\documentclass[final]{siamltex}
\usepackage{epsfig}
\usepackage{graphicx}
\usepackage{graphics}

\begin{document}
\title{Nested  hierarchies in planar graphs}
\author{Won-Min Song\footnotemark[2]\and T. Di Matteo \footnotemark[2]\ \footnotemark[3]
\and Tomaso Aste \footnotemark[2]\ \footnotemark[3]}

\maketitle
\renewcommand{\thefootnote}{\fnsymbol{footnote}}
\footnotetext[2]{Applied Mathematics, Research School of Physical Sciences, The Australian National University, 0200 Canberra, Australia\ ({\tt won-min.song@anu.edu.au}).}
\footnotetext[3]{Department of Mathematics, King's College, The Strand, London, WC2R 2LS, UK.}
\renewcommand{\thefootnote}{\arabic{footnote}}

\begin{abstract}
We construct a partial order relation which acts on the set of 3-cliques of a maximal planar graph $G$ and defines a unique hierarchy. We demonstrate that $G$ is the union of a set of special subgraphs, named `bubbles', that are themselves maximal planar graphs. The graph $G$ is retrieved by connecting these bubbles in a tree structure where neighboring bubbles are joined together by a 3-clique. Bubbles naturally provide the subdivision of $G$ into communities and the tree structure defines the hierarchical relations between these communities.
\end{abstract}



\section{Introduction}
There has been an increasing amount of interest in the study of complex systems via tools of network theory \cite{Boccaletti2006}. Properties such as small-world or scale-free degree distributions have emerged as universal properties of many real complex networks and seemingly these ingredients shape the world we live in \cite{Caldarelli2007,Boccaletti2006}.
In particular, it has been pointed out that the understanding of the organization of local communities is one of the key elements in the study of the structure of complex networks \cite{Girvan2002} and it can shed light on several relevant issues  \cite{Aaron2008}.
One of the underlying assumptions in these studies is that there are local communities and there is a hierarchy among these communities  \cite{Newman2003}.
However, a precise definition of communities and their hierarchy is hard to formulate.\\
In previous works \cite{Tumminello2005,Aste2005,DiMatteo2005,Tumminello2007}, some of the authors proposed a tool for filtering information in complex systems by using planar maximally filtered graphs (PMFG).
This filtering procedure yields to maximal planar graphs\footnote{ A graph is said to be planar if it has an embedding in the plane such that no two edges intersect except at their end points. It is said to be maximal planar if no more edges can be inserted without violating planarity.} that are triangulations of a topological sphere (orientable surface of genus $g=0$) \cite{Diestel2005}.\\
In this paper, we explore ways to characterize the hierarchical structure of maximal planar graphs and we propose a new framework to define communities on these graphs and to extract their hierarchical relation. Planar graphs can display different levels of complexity featuring some of important ingredients of complex networks such as large clustering coefficients, small-world properties and scale-free degree distributions \cite{Andrade2005}. Constituting elements of maximal planar graphs are surface triangles and, more generally, 3-cliques\footnote{A 3-clique is a sub graph made of 3 mutually connected vertices.}. These building blocks also define a class of larger subgraphs, that we name `bubbles', which are themselves maximal planar graphs. We will show in this paper that a hierarchical relationship emerges naturally in planar graphs and it is directly associated to the system of 3-cliques and to the bubble structure.\\
The paper is organized as follows. In {\bf Section \ref{S.1}}, we define a hierarchy among the system of 3-cliques in a maximal planar graph. The concept of hierarchical graph $H$ associated to the  3-cliques hierarchy is introduced in {\bf Section \ref{S.1b}} where it is also shown that $H$ is a forest of rooted trees. This hierarchy is extended to bubbles in {\bf Section \ref{S.2}} where a bubble hierarchical graph $H_b$ is defined and it is shown that $H_b$ is a tree. The generality of the subdivision into bubbles and of the uniqueness of the emerging topology of the bubble hierarchical tree are discussed in {\bf Section \ref{S.2b}}. {\bf Section \ref{S.4}} provides two examples where these hierarchical constructions are applied and discussed in details. Conclusions and perspectives are provided in {\bf Section \ref{S.5}}. In order to improve readability, some of the proofs are reported in {\bf Appendix}.

\begin{figure}[ht!]
\centering
\begin{tabular}{cc}
\epsfig{file=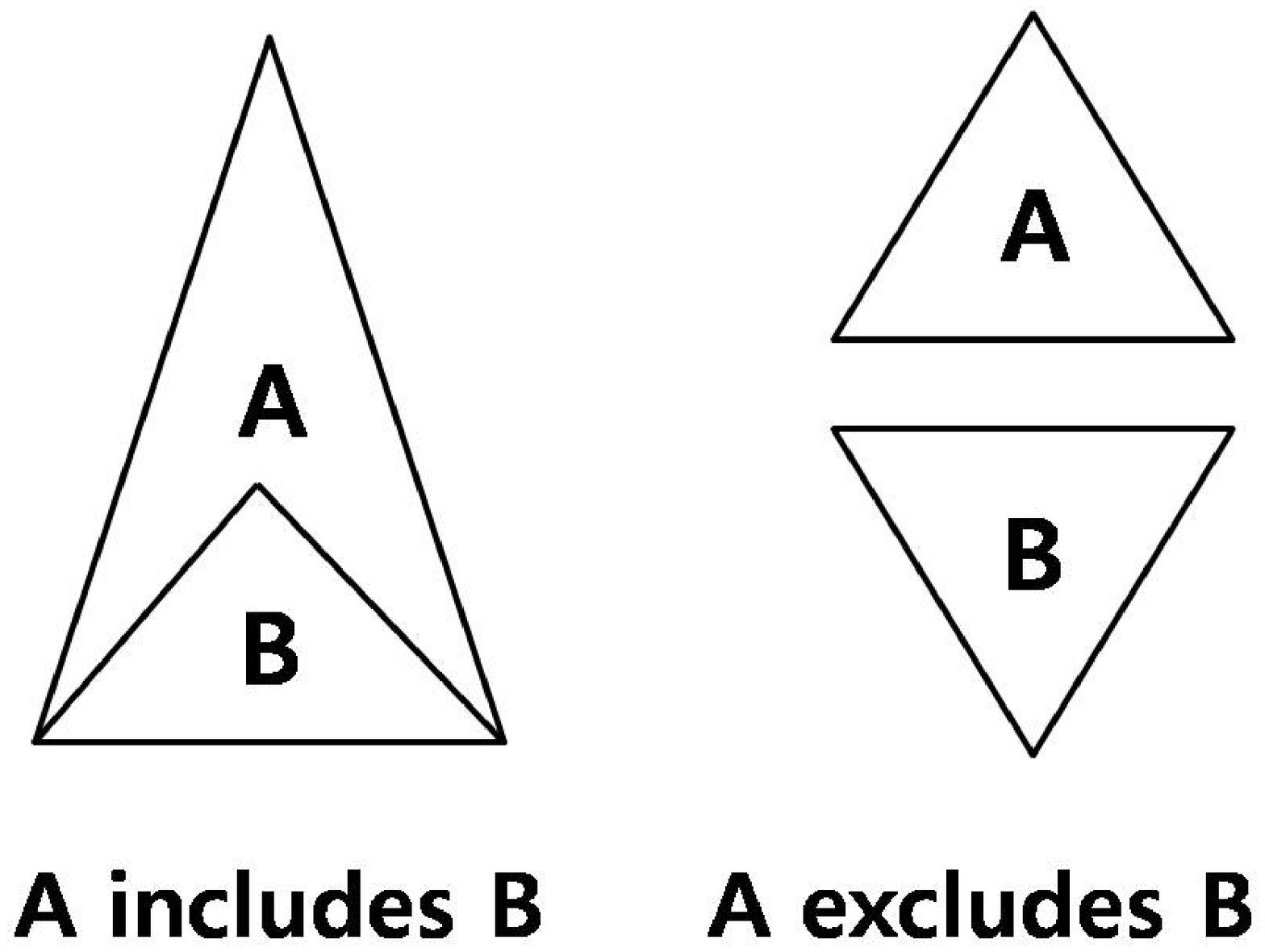,width=5cm,height=4cm}&\epsfig{file=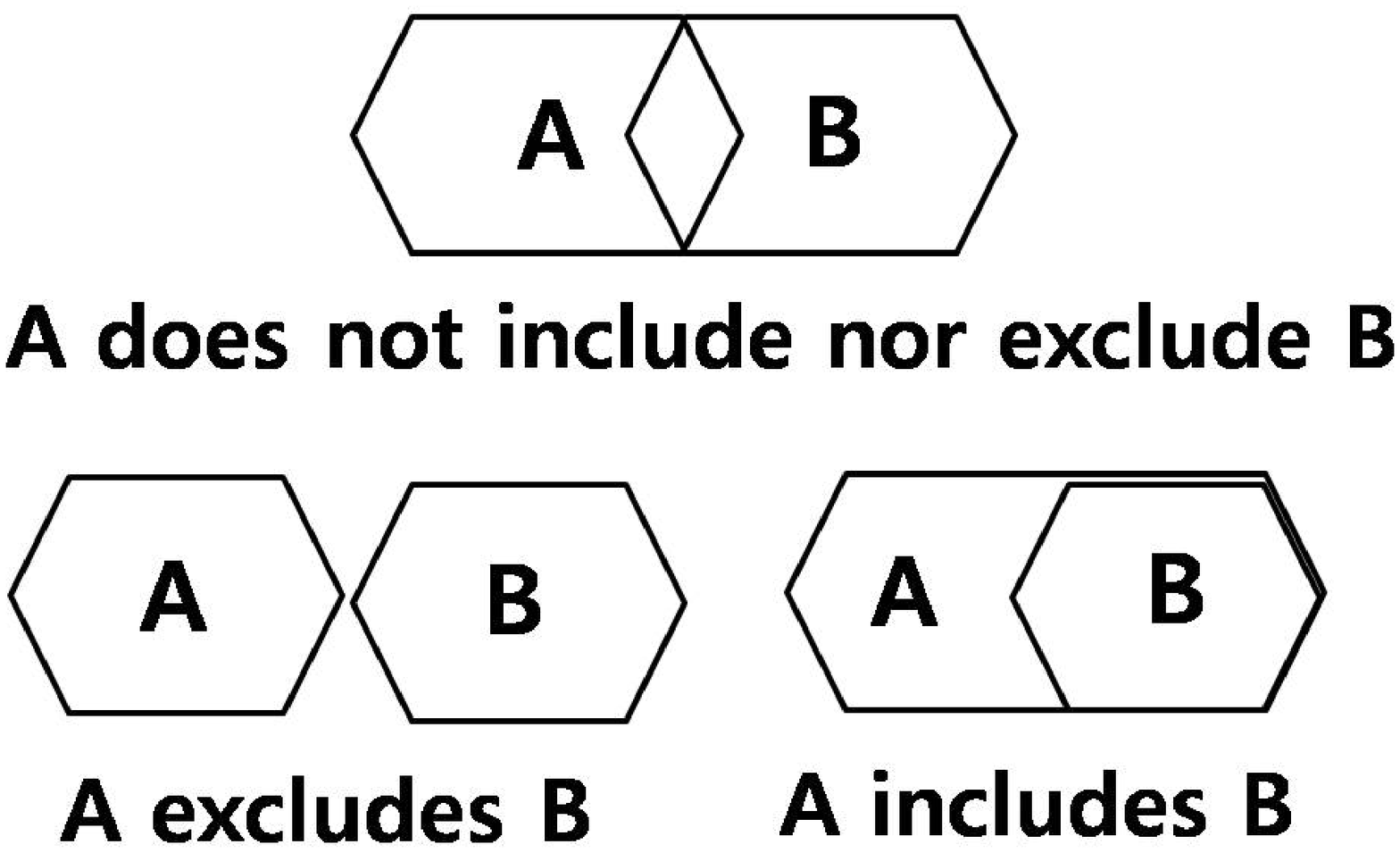,width=5cm,height=4cm}\\
\hbox{(a)}&\hbox{(b)}
\end{tabular}
\caption{(a) Two possible cases for two triangles (3-cliques).
Either one is inclusive of the other or exclusive.
(b) Three possible cases for general cycles larger than a 3-clique showing that there are intermediate cases which are neither inclusive or exclusive.\label{incl}
}
\end{figure}

\section{Hierarchy on Maximal Planar Graphs}\label{S.1}
In this section, we introduce the concept of hierarchy on maximal planar graphs. In particular, we identify a hierarchy as a partially ordered set $(K,\preceq_K)$ where $\preceq_K$ is a binary relation over a set $K$.
We restrict our investigation to maximal planar graphs whose vertices participate to at least one 3-clique and have no multiple edges or self-loops. 3-cliques are the simplest non trivial topological graphs and they are easy to find computationally by looking for the common neighbors of two vertices connected directly by an edge. We will show that the simplicity of 3-cliques provides us with an unique property:
\emph{a 3-clique strictly includes (or excludes) another} (see Fig.~\ref{incl}), a feature that is not true for general cycles. Let us use the notation $G(V,E)$  for a maximal planar graph with  a vertex set $V$ and an edge set $E$. Let us call $K=\{k_1,k_2,...,k_c\}$ the set of 3-cliques in $G(V,E)$, where $k_i$ is the `$i^{th}$' 3-clique. We consider maximal planar graphs which make triangulated surfaces containing at least $3(|V|-2)$ 3-cliques. In order to assign a partial order relation on the set of 3-cliques we proceed in two steps:
\begin{itemize}
\item[]\emph{\bf First} we define the direction of $k_i$ which assigns the \emph{interior} and the \emph{exterior} of the 3-clique;
\item[]\emph{\bf Second} we define the partial order relation $\preceq_K$ from the definition of interior and exterior of $k_i\in K$ on $G$.
\end{itemize}

\subsection{Interior and exterior of a 3-clique}
In a planar graph 3-cliques have a unique property described as follow:
\begin{lemma}[3-clique Removal]\label{cliqueremoval}
Every 3-clique $k_i\in K$ is either a separating or non-separating cycle in $G$.
\end{lemma}\\
Lemma \ref{cliqueremoval} follows from the fact that any cycle is always either a separating or a non-separating cycle in a maximal planar graph \cite{Diestel2005}.
If the vertices $i$ and $j$ are joined by an edge, we shall call the edge $ij$.
Let us recall that a cycle is a graph which consists of a set of distinct vertices $i,j,..,m$ together with a set of distinct edges joining them in cyclic order $ij,...,mi$.
A cycle, as a subgraph of $G$, is said to be \emph{separating} if it divides $G$ into two non-empty subgraphs $S$ and $S'$ that any pair of vertices $i\in S$ and $j\in S'$ is connected by a path which always includes at least one vertex from the cycle.
On the contrary, for a \emph{non-separating cycle}, either of $S$ or $S'$ is empty \cite{Diestel2005}.
\begin{definition}[Interior and Exterior of $k_i$]\label{Orient}
If $S$ and $S'$ are the two subgraphs obtained upon removal of $k_i$, then the one with smaller order\footnote{The order of $S$ is its number of vertices.} is said to be the  \textbf{interior} of $k_i$ and conversely  the \textbf{exterior} of $k_i$ is the one with largest order.
In case the orders of $S$ and $S'$ are equal, the interior is assigned arbitrarily to one of the two subgraphs and exterior to the other.
\end{definition}\\
We denote the interior as $G^i_{in}(V^i_{in},E^i_{in})$ and the exterior as $G^i_{out}(V^i_{out},E^i_{out})$ for each $k_i$.
\begin{definition}[Union and Intersection Operators on Graphs]\label{GraphOper}
Given a graph $G$, let us call two arbitrary subgraphs of $G$ as $S_1$ and $S_2$. Then the union operator $\cup$ on $S_1$ and $S_2$, $S_1\cup S_2$ is defined as follows:
$v\in (S_1\cup S_2)$ if and only if $v\in V_1$ or $v\in V_2$ where $V_1$ and $V_2$ are vertex sets of $S_1$ and $S_2$.\\
Similarly, we define the intersection as follow:
$v\in (S_1\cap S_2)$ if and only if $v\in V_1$ and $v\in V_2$. \end{definition}\\
We say $v_i$ and $v_j$ are connected in $(S_1\cup S_2)$ if and only if $v_iv_j\in E$ and $v_i,v_j\in(V_1\cup V_2)$ where $E$ is the edge set of $G$.\\
We say $v_i$ and $v_j$ are connected in $(S_1\cap S_2)$ if and only if $v_i,v_j\in(V_1\cap V_2)$ and $v_iv_j\in E$.

\begin{theorem}\label{inclusion}
For any two arbitrary 3-cliques $k_i$ and $k_j$, we always have:
\begin{remunerate}
\item If $k_j\subseteq (k_i\cup G^i_{in})$, then $(k_j\cup G^j_{in})\subseteq (k_i\cup G^i_{in})$ or vice versa;
\item otherwise, $G^i_{in}\cap G^j_{in}=\emptyset$.
\end{remunerate}
\end{theorem}
\begin{proof}
Let us call $S^i_{in}=(k_i\cup G^i_{in})$ and $S^j_{in}=(k_j\cup G^j_{in})$ for
convenience.\\
\begin{remunerate}
\item[1.]Firstly, $S^i_{in}$ is maximal planar \cite{Diestel2005}. Therefore, removal of $k_j\subseteq S^i_{in}$ yields to two separate subgraphs $S_1$ and $S_2$ of $S^i_{in}$ so that $S^i_{in}=(S_1\cup k_j\cup S_2)$. Since $k_i\subseteq S^i_{in}$ and $S^i_{in}$ is a disjoint union of $S_1$,$S_2$ and $k_j$, $k_i$ must belong to either of $(S_1\cup k_j)$ or $(S_2\cup k_j)$. Let us arbitrarily choose $k_i\subseteq(k_j\cup S_2)$. Then,
\begin{eqnarray}
G&=(G^i_{in}\cup k_i\cup G^i_{out})&=S^i_{in}\cup G^i_{out}\label{l1}\\
 &  &=(S_1\cup k_j\cup S_2)\cup G^i_{out}\label{l2}\\
 &  &=\left[(S_1\cup k_j)\cup S_2\right]\cup G^i_{out}\label{l3}\\
 &  &=(S_1\cup k_j)\cup(S_2\cup G^i_{out}) \;.
\label{l4}
\end{eqnarray}
We made use of associativity of union operator in the calculation above\footnote{We have omitted the proof for the associativity of the union in Definition~\ref{GraphOper} since it is an immediate consequence of the associativity of union operating on general sets.}. Since $k_i\neq k_j$, there exists at least one vertex of $k_j$ in $S_2$. Then $S_2$ is connected to $G^i_{out}$ via the vertex in $k_i$. Therefore, we can decompose Eq.~\ref{l4} into three disjoint union of connected subgraphs $S_1$,$k_j$ and $(S_2\cup G^i_{out})$ as:
\begin{eqnarray}
G&=S_1\cup k_j\cup(S_2\cup G^i_{out})\;\;.
\label{l5}
\end{eqnarray}
In Eq.\ref{l5}, it is immediate that $S_1$ and $(S_2\cup G^i_{out})$ are the disjoint subgraphs of $G$ realized by removal of $k_j$ as in Lemma~\ref{cliqueremoval}. Comparing the orders, $|S_1|\leq|G^i_{in}|\leq|G^i_{out}|\leq|(S_2\cup G^i_{out})|\Rightarrow|S_1|\leq|(S_2\cup G^i_{out})|$. By Definition~\ref{Orient}, it is immediate that $G^j_{in}=S_1$ and $G^j_{out}=(S_2\cup G^i_{out})$.\\
\item[2.] Suppose $(G^i_{in}\cap G^j_{in})\neq\emptyset\Rightarrow \exists v_o\in(G^i_{in}\cap G^j_{in})$. And suppose $k_j\not\subseteq (k_i\cup G^i_{in})$ and vice versa so that $\exists v_j\in k_j\hbox{ such that }v_j\in G^i_{out}$. Then $v_o$ and $v_j$ are connected via $k_i$. However, $v_j$ is connected to any vertices in $G^j_{in}$ including $v_o$ without $k_i$ since we have assumed $k_i\not\subseteq(k_j\cup G^j_{in})$. Therefore, by contradiction, $G^i_{in}\cap G^i_{in}=\emptyset$.
\end{remunerate}
\end{proof}
\begin{corollary}\label{InCliq}
Given a 3-clique $k_i$ and a vertex $v\in G^i_{in}$, all 3-cliques $k_j$ in which $v\in k_j$ satisfy $k_j\subseteq(k_i\cup G^i_{in})$.
\end{corollary}
\begin{proof}
By Lemma~\ref{cliqueremoval}, $v$ is connected to $G^i_{out}$ via $k_i$. This implies that any vertex in $G^i_{in}$ is connected to $G^i_{out}$ by paths made of at least 2 edges. Then, any 3-clique involving $v$ cannot have vertices from $G^i_{out}$ since vertices in a 3-clique are connected by single edges, not paths of greater lengths. Therefore, $k_j\subseteq(k_i\cup G^i_{in})$.
\end{proof}
\subsection{Partial Order Relation on the set of 3-cliques}
\begin{figure}[ht!]
\centering
\begin{tabular}{cc}
\epsfig{file=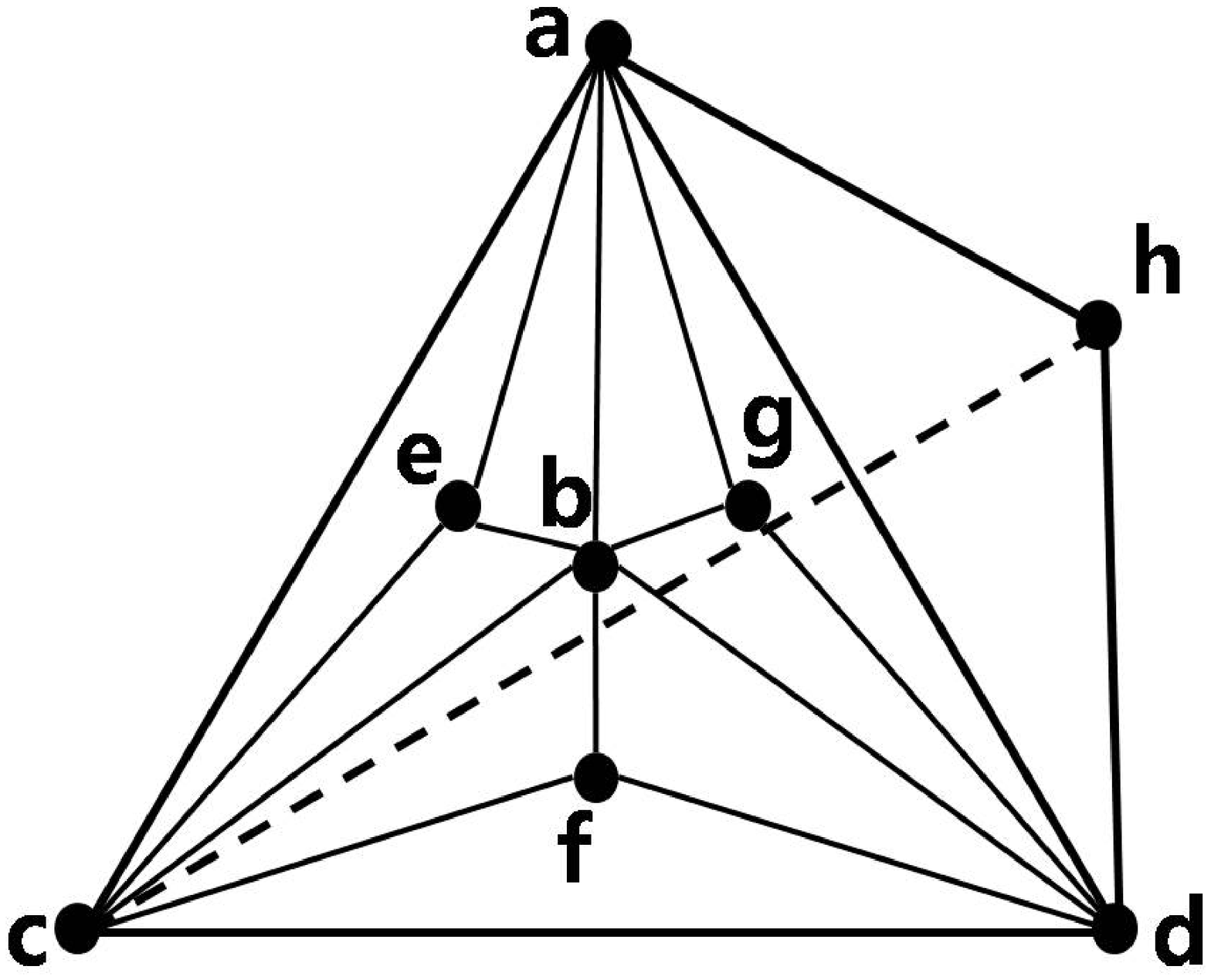,height=4.5cm,width=5.5cm}&\epsfig{file=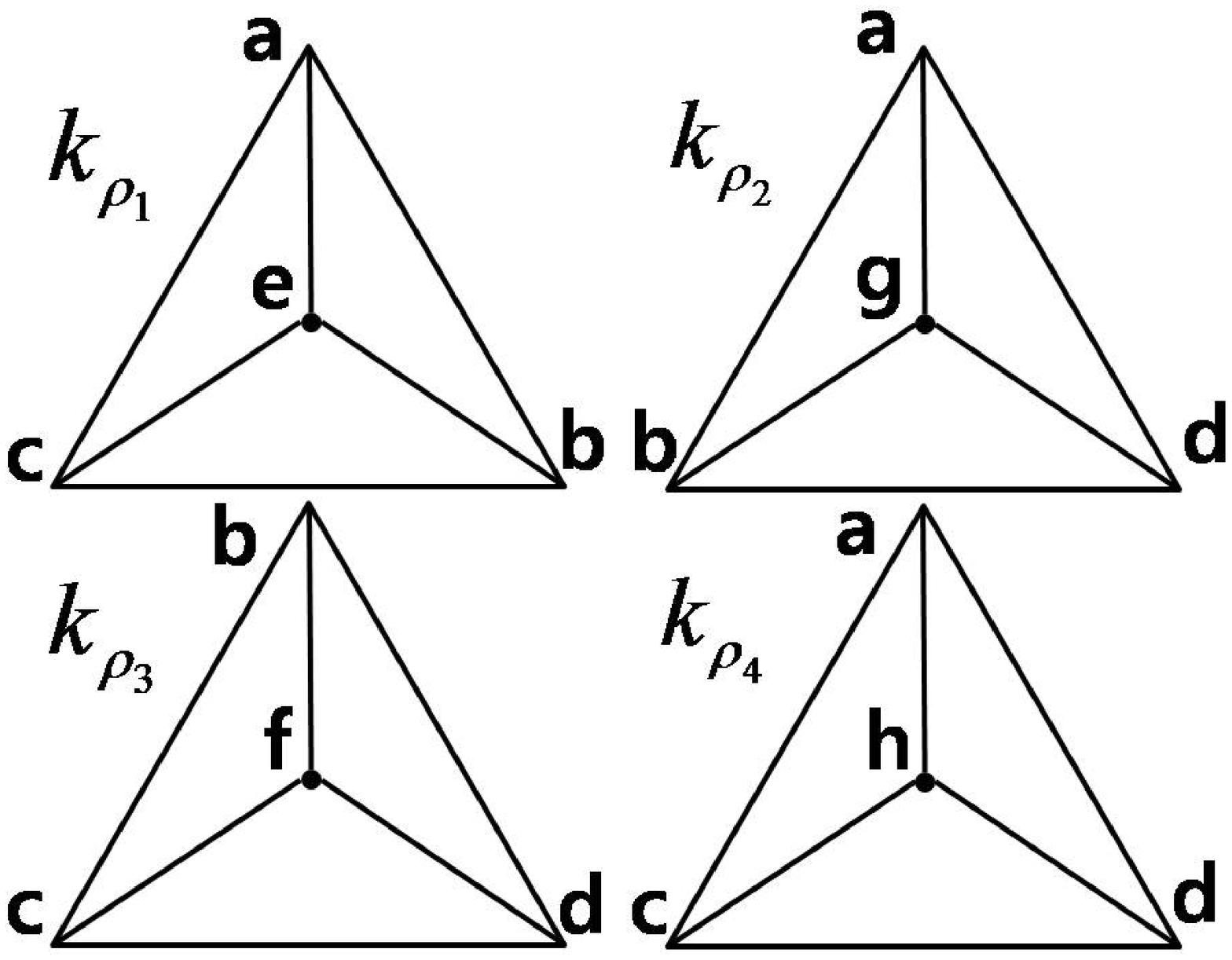,height=4.5cm,width=5.5cm}\\
\hbox{(a)}&\hbox{(b)}
\end{tabular}
\caption{(a) A maximal planar graph which has 4 maximal elements. (b) We have dissembled the maximal planar graph into the four maximal 3-cliques. The 3-cliques are (labels in the brackets denote the vertices of corresponding 3-clique): $k_{\rho_1}=(a,b,c)$ ,$k_{\rho_2}=(a,b,d)$, $k_{\rho_3}=(b,c,d)$, $k_{\rho_4}=(a,c,d)$. All of these 3-cliques have the same size of interior, which is one vertex as shown.\label{maximaleg}}
\end{figure}

Let us here introduce the relation $\preceq_K$ for the 3-cliques in $K$.
\begin{definition}[Relation between two 3-cliques]\label{Order}
$k_i\preceq_K k_j$ if and only if $(k_i\cup G^i_{in})\subseteq (k_j\cup G^j_{in})$.
\end{definition}
We now show with the following theorem that $\preceq_K$  is a partial order relation.
\begin{theorem}[Partial order relation between 3-cliques]\label{axioms}
The partially ordered set $(K,\preceq_K)$ satisfies the three axioms:
\begin{remunerate}
\item Reflexivity: $k_i\preceq_K k_i$ for $\forall k_i\in K$,
\item Antisymmetry: $k_i\preceq_K k_j$ and $k_j\preceq_K k_i$ implies $k_i=k_j$,
\item Transitivity: $k_i\preceq_K k_j$ and $k_j\preceq_K k_m$ implies $k_i\preceq_K k_m$.
\end{remunerate}
\end{theorem}
\begin{proof}
\begin{enumerate}
\item[Reflexivity] For all $i$, $(k_i\cup G^i_{in})\subseteq (k_i\cup G^i_{in})$. Therefore reflexivity holds.
\item[Antisymmetry] If $k_i\preceq_K k_j$, then $(k_i\cup G^i_{in})\subseteq (k_j\cup G^j_{in})$ by definition. Similarly, if $k_j\preceq_K k_i$, then $(k_j\cup G^j_{in})\subseteq (k_i\cup G^i_{in})$. These imply $(k_i\cup G^i_{in})=(k_j\cup G^j_{in})$. However, this does not mean $k_i=k_j$.\\
    So, suppose that $k_i\neq k_j$ to prove by contradiction. Then $G^i_{in}\neq G^j_{in}$ and $G^i_{out}=G^j_{out}$. This implies that there exists a vertex $v_i\in k_i$ such that $v_i\in G^j_{in}$.
 By Lemma~\ref{cliqueremoval}, $v_i$ is also directly connected to $G^i_{out}$. Then, $G^j_{out}$ is directly connected to $G^j_{in}$ without  $k_j$ since $G^i_{out}=G^j_{out}$. This violates the planarity of $G$, therefore $k_i=k_j$.
\item[Transitivity] If $k_i\preceq_K k_j$ and $k_j\preceq_K k_m$, then $(k_i\cup G^i_{in})\subseteq (k_j\cup G^j_{in})$ and $(k_j\cup G^j_{in})\subseteq (k_m\cup G^m_{in})$. Therefore, $(k_i\cup G^i_{in})\subseteq (k_m\cup G^m_{in})$. So $k_i\preceq_K k_m$ is true.
\end{enumerate}
\end{proof}
Note that  $(K,\preceq_K)$ is a partially ordered set (a poset) and therefore, differently from a total order, there are some elements in $K$ that might not be related to each other through $\preceq_K$. It is known that, for any finite poset, one can find a set of maximal elements which are not smaller that any other element in the set \cite{Jech2003}.
For the set of 3-cliques in $G$, we have therefore the following theorem:
\begin{theorem}[Maximal Elements]\label{maximal}
$(K,\preceq_K)$ always has at least one maximal element, and can have more than one maximal element.
\end{theorem}
\begin{proof}
If a poset has a finite order, then it always has at least one maximal element \cite{Jech2003}.
Since $K$ is a finite set, there is at least one maximal elements.
In order to prove that there can be many maximal elements, we provide an example in Fig.~\ref{maximaleg} that has several maximal elements.
\end{proof}

\section{Hierarchical Graph for the 3-cliques}\label{S.1b}
We can now associate to the hierarchical relation between 3-cliques a graph where the vertices are the 3-cliques and the \emph{directed} edges connect 3-cliques accordingly to the poset relation $\preceq_K$. Let us first formalize closest elements in the poset by defining covering elements, so that we can associate the elements by edges.
\begin{definition}[Covering elements]
Given the poset $(K,\preceq_K)$, an element  $x\in K$ is said to cover $y\in K$ if $y\preceq_K x$ and there is no other element $z\in K$ such that $y\preceq_K z$ and $z\preceq_K x$ \cite{Jech2003}.
\end{definition}
\begin{definition}[Hierarchical Graph]\label{H}
The hierarchical graph $H(K,E_k)$, has vertex set $K$ and a directed edge $\overrightarrow{k_jk_i}\in E_k$ from $k_j$ to $k_i$ is a pair such that $k_i$ covers $k_j$ in $(K,\preceq_K)$. We shall call $\overrightarrow{k_jk_i}$ outgoing edge from $k_j$ and incoming edge to $k_i$.
\end{definition}
Having defined the hierarchical graph, one can define adjacent elements as {\bf neighbors} which is a general term in graph theory \cite{Diestel2005}, and further characterize the properties in the language of graph theory.
\begin{definition}[Neighbors]\label{Neigh}
Given a vertex $v_o$ in an undirected graph $X$, we say $v_1,v_2,..$ are neighbors of $v_o$ if they share edges with $v_o$.
If $X$ is a directed graph with incoming and outgoing edges  then we call $v_1,v_2,..$ (incoming/outgoing) neighbors if they share an (incoming/outgoing) edge at $v_o$.
\end{definition}
\begin{theorem}[Incoming Neighbors]\label{SeparateCliq}
For an arbitrary 3-clique $k_i$ with incoming neighbors $k_j,k_l,..$ in $H$, one has $G^j_{in}\cap G^l_{in}=\emptyset$.
\end{theorem}
\begin{proof}
Suppose $k_j\subseteq(k_l\cup G^l_{in})$ or $k_l\subseteq(k_j\cup G^j_{in})$. This implies that it is either $k_j\preceq_K k_l$ or $k_l\preceq_K k_j$, therefore $k_i$ does not cover $k_l$ or $k_j$. By Definition~\ref{H}, this is against the definition of the hierarchical tree. By contradiction, this implies $k_j\not\subseteq(k_l\cup G^l_{in})$ and $k_l\not\subseteq(k_j\cup G^j_{in})$. Therefore, by Theorem~\ref{inclusion}, $G^j_{in}\cap G^l_{in}=\emptyset$.
\end{proof}
\begin{theorem}\label{Hedge}
Any vertex in $H$ can have (a) several incoming edges, but (b) no more than one outgoing edge.
\end{theorem}
\begin{figure}[ht!]
\centering
\begin{tabular}{cc}
\epsfig{file=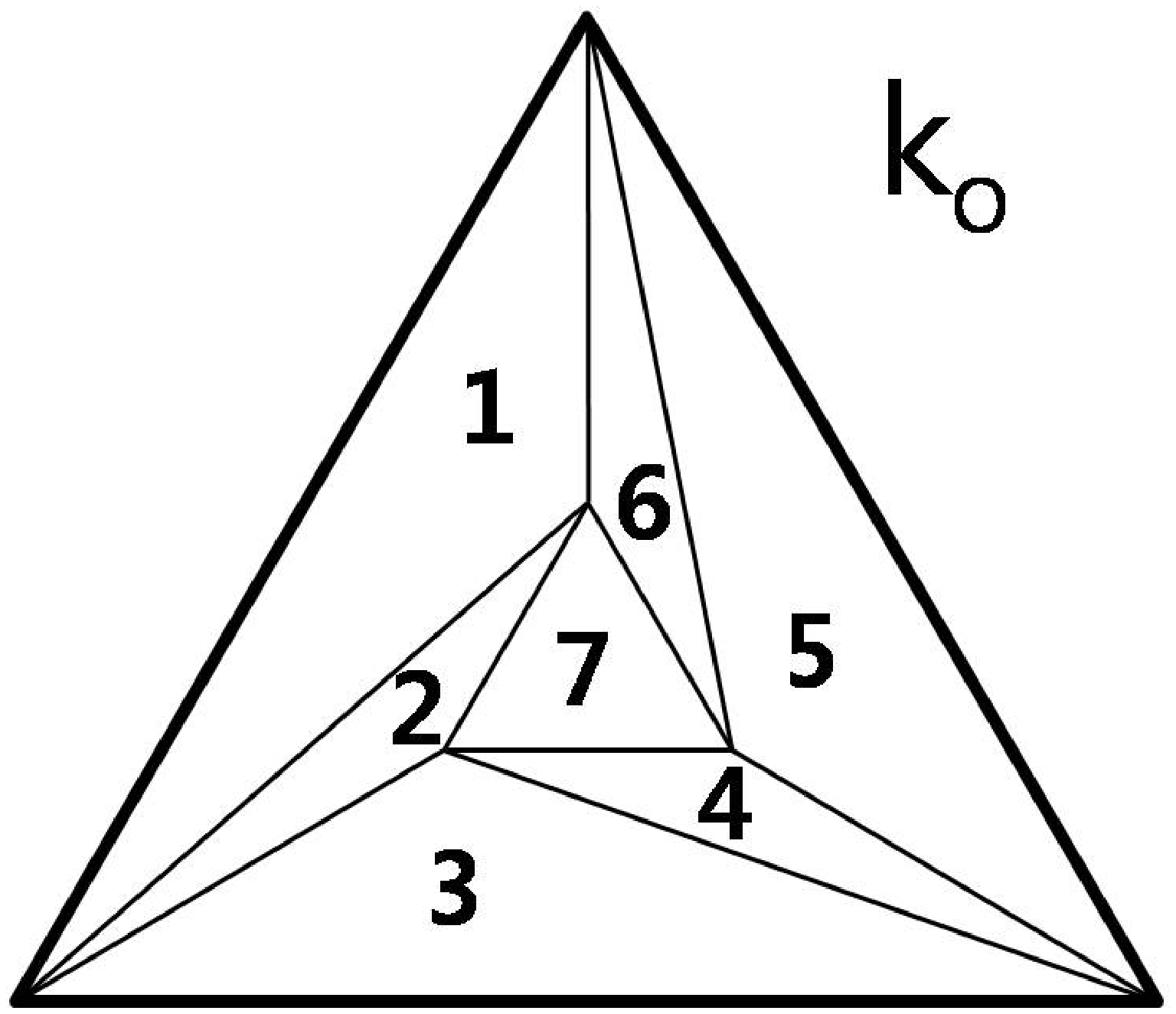,height=4cm,width=5cm}&\epsfig{file=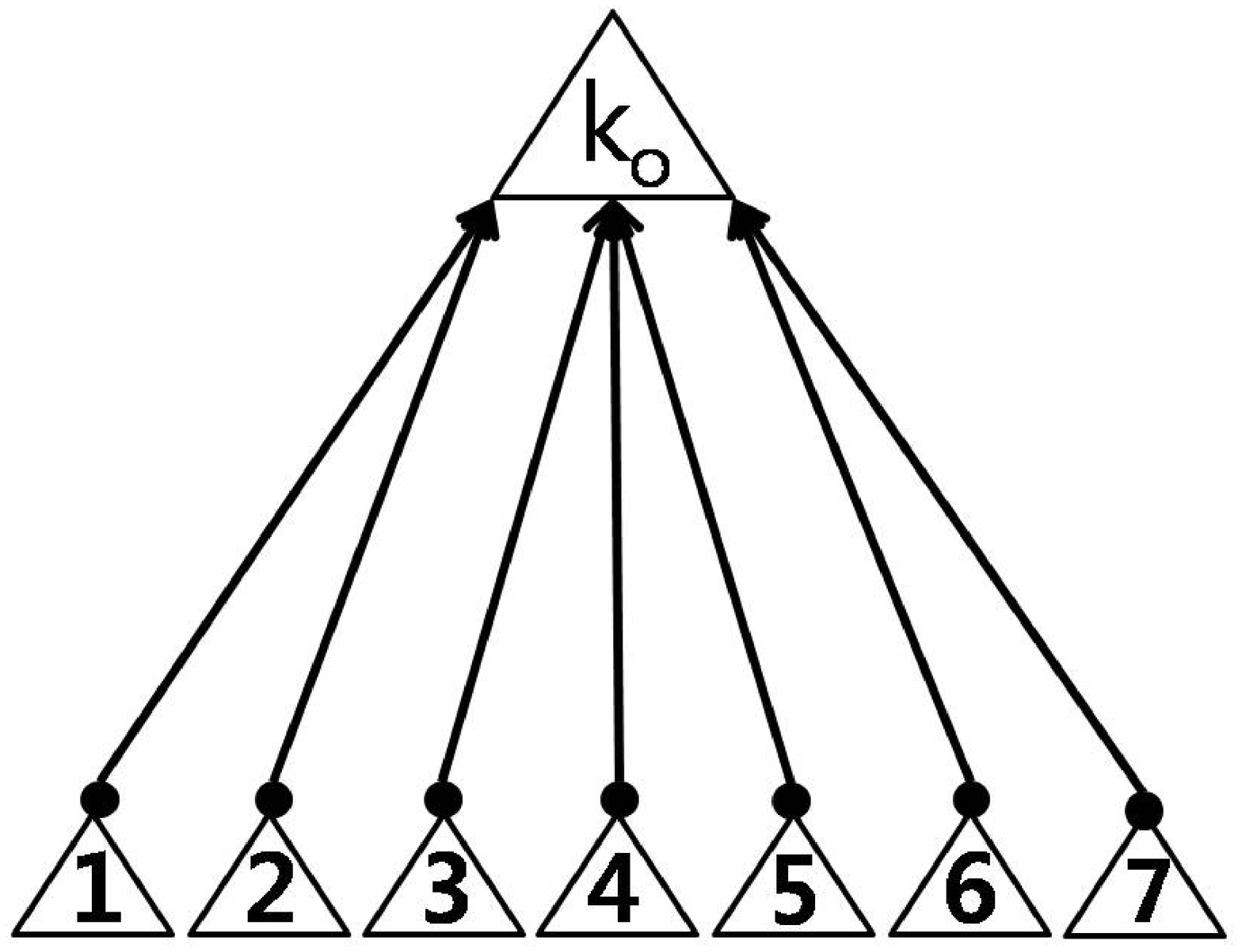,height=4cm,width=5cm}\\
\hbox{(a)}&\hbox{(b)}
\end{tabular}
\caption{(a) The interior of a 3-clique $k_o$ can be triangulated by non-separating triangles which have empty interiors. (b) Each of them form a different incoming edge at the vertex $k_o$  in the hierarchical graph $H$.\label{Several}}
\end{figure}
\begin{proof}
\begin{enumerate}
\item[(a)]In other words, we need to prove that there can be a number of 3-cliques which reside in the interior of an arbitrary 3-clique $k_o$.
Let us present an example. Given $k_o$, let us suppose that it's interior corresponds to the `inside of the triangle' in Fig.~\ref{Several}(a).
In Fig.~\ref{Several}(a) there are seven non-separating 3-cliques in the interior of $k_o$, these 3-cliques are incoming edges in $H$ by Definition~\ref{H} as shown in (b).
\item[(b)]
Let us first show that there can be one outgoing neighbor in $H$. Fig.~\ref{Several} is an example that, indeed, $k_o$ is the only outgoing neighbor of the  3-cliques $1,2,3,...,7$.\\
Now, proving by contradiction, let us suppose that there are more than one outgoing edges, and say there are two outgoing edges without loss of generality.
Let us call $k_1$ and $k_2$ the outgoing neighbors of $k_o$. Then $[k_o\subseteq(k_1\cup G^1_{in})]\wedge[k_o\subseteq(k_2\cup G^2_{in})]\Rightarrow G^1_{in}\cap G^2_{in}\neq\emptyset$ since $(k_o\neq k_1)\wedge(k_o\neq k_2)$.\\
    We have two possible cases between $k_1$ and $k_2$ as suggested by Corollary~\ref{inclusion}:
    \begin{remunerate}
    \item[(i)]$k_1\subseteq(k_2\cup G^2_{in})\Rightarrow k_1\preceq_K k_2$. Since $k_o\preceq k_1$ by the assumption, $k_o$ is not covered by $k_2$. This violates the assumption that $k_2$ covers $k_o$. The same argument holds when $k_2\subseteq(k_1\cup G^1_{in})$.
    \item[(ii)]$[k_1\not\subseteq(k_2\cup G^2_{in})]\wedge[k_2\not\subseteq(k_1\cup G^1_{in})]\Rightarrow G^1_{in}\cap G^2_{in}=\emptyset$. But, this also violates the initial assumption that $G^1_{in}\cap G^2_{in}\neq\emptyset$.
    \end{remunerate}
    Therefore, there cannot be two outgoing neighbors since the assumption of two outgoing neighbors yields contradiction. The same argument holds for general cases of many outgoing neighbors $k_1,k_2,...$. Hence there cannot be more than one outgoing neighbor in $H$ for all 3-cliques.
 \end{enumerate}
\end{proof}
\begin{corollary}[Forest of Rooted Hierarchical Trees]\label{HForest}
$H$ is a forest of rooted trees where the number of trees corresponds to the number of maximal elements in $(K,\preceq_K)$ and each maximal element is the root of a tree.
\end{corollary}
\begin{proof}
\begin{enumerate}
\item[No Cycle]
In order to be a forest or a tree, $H$ must have no cycles. In order to prove that $H$ does not possess any cycle, let us suppose that there exists a cycle which is made of a set of distinct edges. Without loss of generality, let us suppose that the cycle is of order 3 and expressed as $\overrightarrow{k_ik_j}\overrightarrow{k_jk_m}\overrightarrow{k_mk_i}$. By definition of the edges in $H$, this implies $k_i\preceq_K k_j$, $k_j\preceq_K k_m$. By transitivity of $(K,\preceq_K)$, this implies $k_i\preceq_K k_m$. However, the cycle also implies $k_m\preceq_K k_i$. By reflexivity of $(K,\preceq_K)$, this implies $k_i=k_m$. This is against the assumption that the cycle is of order 3. Therefore, the assumption is false. The same argument using the transitivity and reflexivity applies to any cycle of order grater than 3. Therefore there does not exist any cycle in $H$.
\item[Forest] In order to prove that $H$ can be a forest of many trees, it is sufficient to show that $H$ can be disconnected. This implies that there exist two 3-cliques $k_i$ and $k_j$ which do not have a connecting path in $H$. We will use the maximal elements in $(K,\preceq_K)$ to prove the disconnectedness of $H$. Let us suppose that $(K,\preceq_K)$ possesses more than one maximal element (which is a possible case by Theorem~\ref{maximal}), and let there be two maximal elements $k_{\rho_i}$ and $k_{\rho_j}$. Suppose they are connected by a path. By transitivity, this implies $k_{\rho_i}\preceq_K k_{\rho_j}$ or vice versa. This is false because these are maximal elements. Therefore all maximal elements of $(K,\preceq_K)$ are disconnected. This implies that $H$ is a forest made of a number of tree greater or equal than the number of maximal elements.
\item[Roots]
In $H$, a root is a 3-clique which does not possess outgoing edges. Clearly, the maximal element does not possess any outgoing edges by definition. Unless a 3-clique $k_i$ is a maximal element, then it is not a root since there exists always another 3-clique $k_j$ which is $k_i\preceq_K k_j$. This also proves that the number of trees is equal to to the number of maximal elements.
\end{enumerate}
\end{proof}
\begin{definition}[Nested Community]
A tree in $H$ is a nested community.
\end{definition}
\begin{definition}[Nesting Depth]
The path length between a 3-clique $k_i$ and its corresponding root in a nested community is the nesting depth of $k_i$.
\end{definition}\\
We denote the nesting depth as $h(k_i)$. The graph $H$ provides a valuable instrument to study hierarchy in maximal planar graph  $G$. Specifically we have two important properties that can be used to classify hierarchies in  $G$:
{\bf first}, we have a natural division of $G$ into a system of subgraphs associated with the rooted trees;
{\bf second}, the nesting depth provides us with an instrument to  further distinguish between the various branches inside the trees as a function of increasing topological distance from the root. On the other hand, by construction, the structure of $H$ and the hierarchy depend on the definition of interior and exterior for the separating cliques. Although, well defined, such a dependence on the clique direction can be source of artificial hierarchical positioning of the cliques in $H$. In the next section, we show that we can also eliminate the dependence on clique direction by extending the hierarchy to bubbles.
\section{Hierarchy on bubbles: extension of 3-clique hierarchy}\label{S.2}
Here we extend the concept of nested hierarchy from 3-cliques to larger portions of planar graph that we shall call `bubbles'. This extension has the advantage to produce a connected hierarchical graph with a topological structure that is independent on the choice of interior for the 3-clique.
\begin{definition}[Imaginary 3-clique]\label{ImagCliq}
We define an imaginary 3-clique $k_{imag}$ whose interior $G^{imag}_{in}$ is $G$.
\end{definition}
Let us denote $K'=K\cup\{k_{imag}\}$.
\begin{corollary}\label{CliqHier}
If we define $\preceq_{K'}$ as follows:
\begin{center}
For $k_i,k_j\in K'$, $k_i\preceq_{K'}k_j$ if $(k_i\cup G^i_{in})\subseteq(k_j\cup G^j_{in})$.
\end{center}
Then $(K',\preceq_{K'})$ is a partially ordered set with a single maximal element $k_{imag}$.
\end{corollary}
\begin{proof}
The order relation $\preceq_{K'}$ is identically defined to $\preceq_K$ except that it acts on an extended set $K'$. Thus, the proof to show that $(K',\preceq_{K'})$ is a poset is trivial to that of $(K,\preceq_{K})$ except that we need to show that $k_{imag}$ satisfies the axioms for a poset.\\
The reflexivity holds since $(k_{imag}\cup G^{imag}_{in})\subseteq(k_{imag}\cup G^{imag}_{in})$. The antisymmetry holds since $(k_{imag}\cup G^{imag}_{in})$ is the only graph which can have itself as a subgraph. The transitivity holds since it follows naturally from the transitivity of the operator $\subseteq$.\\
It is immediate from Definition~\ref{ImagCliq} that $k_i\preceq_K k_{imag}$ and $k_i\not\preceq_K k_{imag}$ for all $k_i\in K$, therefore $k_{imag}$ is the only maximal element in $(K',\preceq_{K'})$.
\end{proof}
\subsection{Bubble Hierarchy}
Let us begin by formally define a `bubble':
\begin{definition}[Bubble]\label{bubble}
A bubble $b$ is a maximal planar graph whose 3-cliques are non-separating cycles.
\end{definition}\\
A bubble is a special class of maximal planar graphs where the set of all the 3-cliques are all triangular faces. This implies that each 3-clique is a maximal element as well as a minimal element. A bubble has therefore the simplest hierarchical structure.
Hereafter we use the concept of bubble in order to analyze $G$  as made of a set of bubbles joined by separating 3-cliques. To this end, we search for bubbles in $G$ by making use of the property that each 3-clique is maximal as well as minimal. We can also define a hierarchy for them by making use of $(K,\preceq_K)$ once we describe $G$ in terms of bubbles.
\begin{theorem}[Bubbles in $G$]\label{BubMake}
Given a 3-clique $k_i\in K$ which has incoming edges in $H$, its graph union with the neighbor 3-cliques $k_j,k_m,..$  is a bubble $b_i$.
\end{theorem}
\begin{proof}
The proof is given in Appendix \ref{A.0}.
\end{proof}
Note that Theorem~\ref{BubMake} \emph{does not} realize bubbles made by only one 3-clique since it always takes graph union of more than one 3-clique.
Let us call $k_i$ as \emph{root 3-clique of $b_i$} and the set of all bubbles obtained by Theorem~\ref{BubMake} as $B$.
\begin{theorem}[Maximal 3-cliques]\label{MaxBub}
The graph union of maximal 3-cliques of $(K,\preceq_K)$ is a bubble.
\end{theorem}
\begin{proof}
The proof is given in appendix  \ref{A.1}.
\end{proof}
Let us call this bubble made of the maximal 3-cliques as $b_{\rho}$.
\begin{corollary}[All bubbles]\label{AllBub}
The union $B'=B\cup\{b_{\rho}\}$ is the set of {all bubbles} we can find in $G$.
\end{corollary}
\begin{proof}
The proof is given in Appendix \ref{A.2}.
\end{proof}
This is a very important result because it clarifies that the `bubbles' are defined independently from the 3-clique hierarchy.
\begin{definition}[Bubble Hierarchy]\label{BubHier}
Let $B'$ be a set of bubbles realized by Theorems~\ref{BubMake} and \ref{MaxBub}.
Then let us define the relation $\preceq_{B'}$ between $b_i$ and $b_j$ as follows:
\begin{center}
$b_i\preceq_{B'} b_j$ if $k_i\preceq_{K'} k_j$.
\end{center}
\end{definition}
\begin{theorem}\label{BubPoset}
$(B',\preceq_{B'})$ is a partially ordered set with a single maximal element $b_{\rho}$.
\end{theorem}
\begin{proof}
The proof is trivial from that of $(K',\preceq_{K'})$ because of the one-to-one correspondence between $B'$ and $K'$.
\end{proof}
\section{Bubble Hierarchical Tree}\label{S.2b}
\begin{definition}[Bubble Hierarchical Tree]\label{BubTree}
If $b_i$ covers $b_j$ in $(B',\preceq_{B'})$, then they are connected in the hierarchical tree $H_b(B',E_b)$ by a directed edge $\overrightarrow{b_ib_j}$.
\end{definition}
\begin{corollary}\label{RootTree}
$H_b$ is a single rooted tree whose root is $b_{\rho}$.
\end{corollary}
\begin{proof}
We simply illustrate the proof for Corollary~\ref{RootTree} since it resembles the proof for Corollary~\ref{HForest}.
Indeed, $H_b$ is a forest because does not possess any loop since loops disrupt the planarity.
$H_b$ is a single rooted tree because the relation $b_i\preceq_{B'}b_{\rho}$ holds for all $b_i\in B$.
\end{proof}
As suggested from Corollary~\ref{AllBub}, the bubble hierarchical tree $H_b(B',E_b)$ has the very important property that the connection topology (non-directed edges) is independent on the definition of interior/exterior  of the 3-cliques. Indeed, given a separating 3-clique $k_i$, by definition, it divides $G$ into two maximally planar subgraphs  $(k_i\cup G^i_{in})$ and $(k_i\cup G^i_{out})$. From Theorem~\ref{BubMake}, a 3-clique in a bubble is either: (i) an incoming neighbor at some 3-clique in $H$, or (ii) an outgoing neighbor at other 3-cliques in $H$. This implies that any separating 3-clique belongs to two bubbles, hence the topology of connection between bubbles does not depend on the definition of interior but depends on whether two bubbles share a common 3-clique. Let us now formalize this property more precisely with the two following corollaries.
\begin{corollary}[Separating 3-cliques in $H$]\label{CliqSep}
A 3-clique is separating in $G$ if and only if it has a non-empty set of incoming neighbors in $H$.
\end{corollary}
\begin{proof}
The proof consists of two parts: proving in the forward and backward directions. (i) A 3-clique $k_i$ is separating if it has incoming neighbors in $H$, and (ii) $k_i$ has incoming neighbors in $H$ if it is a separating 3-clique in $G$. \\
{\bf  (i)} Suppose $k_i$ has incoming neighbors in $H$. This implies that $k_i$ has a non-empty interior. By the definition of interior and exterior in Definition~\ref{Orient}, $k_i$ has non-empty exterior as well since it must have a larger subgraph of $G$, $G^i_{out}$, than $G^i_{in}$. Therefore $k_i$ is separating.\\
{\bf  (ii)} Suppose $k_i$ is a separating 3-clique. Then $k_i$ has non-empty interior $G^i_{in}$. Since $G$ is maximally planar, any vertex in $G^i_{in}$ belongs to at least one 3-clique that are not $k_i$. By Corollary~\ref{InCliq}, this 3-clique belongs to $(k_i\cup G^i_{in})$. Therefore, $k_i$ has incoming neighbors in $H$.\\
\end{proof}
Let us now use the above property to state that two adjacent bubbles in $H_b$ always share a common 3-clique. 
\begin{corollary}\label{CommonCliq}
$b_i$ is an incoming or outgoing neighbor at $b_j$ in $H_b$ if and only if the two bubbles $b_i$ and $b_j$ share a common 3-clique.
\end{corollary}
\begin{proof}
Let us proceed in two steps:
\begin{remunerate}
\item[(i)] If $b_i$ is an incoming or outgoing neighbor at $b_j$ in $H_b$, then they share one and only one common 3-clique;
\item[(ii)] If $b_i$ and $b_j$ share a common 3-clique, then $b_i$ is an incoming or outgoing neighbor at $b_j$.
\end{remunerate}
{\bf First}, let us suppose that $b_i$ is an incoming neighbor at $b_j$.
By the definition of edges in $H_b$ (Definition~\ref{BubTree}), there are two 3-cliques $k_i\subset b_i$ and $k_j\subset b_j$ that $k_i$ is incoming neighbor at $k_j$ in $H$. By Theorem~\ref{BubMake}, $k_i\subset b_j$, hence they share at least one 3-clique. However, they cannot share more than one 3-clique. Suppose they do. And let us call these two 3-cliques as $k_i$ and $k'_i$. Since $b_i\preceq_{B'}b_j$, $k_i,k'_i\preceq_K k_j$, $k_i$ and $k'_i$ are not maximal. By Corollary~\ref{AllBub}, $b_i$ is then not a union of maximal 3-cliques but union of $k_i$ and its incoming neighbors in $H$. Therefore $k'_i\preceq_K k_i$, hence $k'_i$ is not covered by $k_j$. Then, by Corollary~\ref{AllBub}, $k'_i$ cannot be a subgraph of $b_j$. This violates the assumption that $k'_i$ is the subgraph of $b_j$.
Therefore, there is always one and only one common 3-clique between the neighboring bubbles $b_i$ and $b_j$ in $H_b$.\\
{\bf Second}, let us suppose that two bubbles $b_i$ and $b_j$ share a common 3-clique and call it $k$.
Suppose that $b_i$ is not an incoming neighbor in $H_b$. Then there exists a bubble $b_l$ such that $(b_l\neq b_i)\wedge(b_l\neq b_j)$, and $(b_i\preceq_{B'}b_l)\wedge(b_l\preceq_{B'}b_j)$. Since $k$ is a 3-clique of $b_i$, it satisfies $k\preceq_K k_i$. Similarly, $k_i\preceq_K k_l$ since $b_i\preceq_{B'}b_l$. Since $k\neq k_i$, $k$ is not covered by $k_j$, therefore $k$ is not a 3-clique of $b_j$. This violates our assumption that $k$ is the common 3-clique of $b_i$ and $b_j$, therefore $b_j$ covers $b_i$, hence they are connected in $H_b$.
\end{proof}
We have therefore proved that neighboring bubbles in $H_b$ are always connected through one and only one common 3-clique and vice versa any separating 3-clique in $H$ is always shared by two bubbles.
\begin{figure}[t]
\centering
\begin{tabular}{ccc}
\epsfig{file=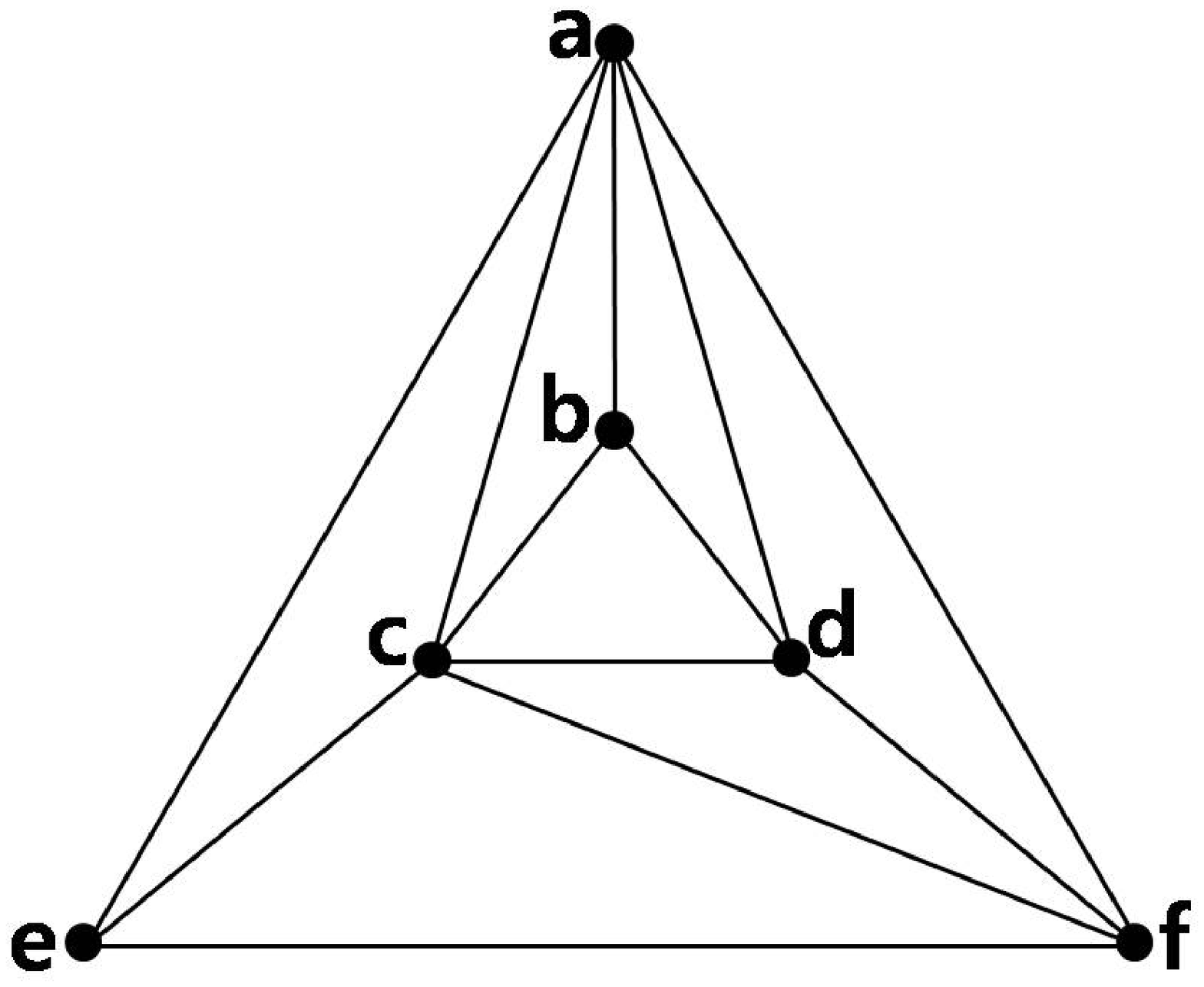,height=3cm,width=4cm}&\epsfig{file=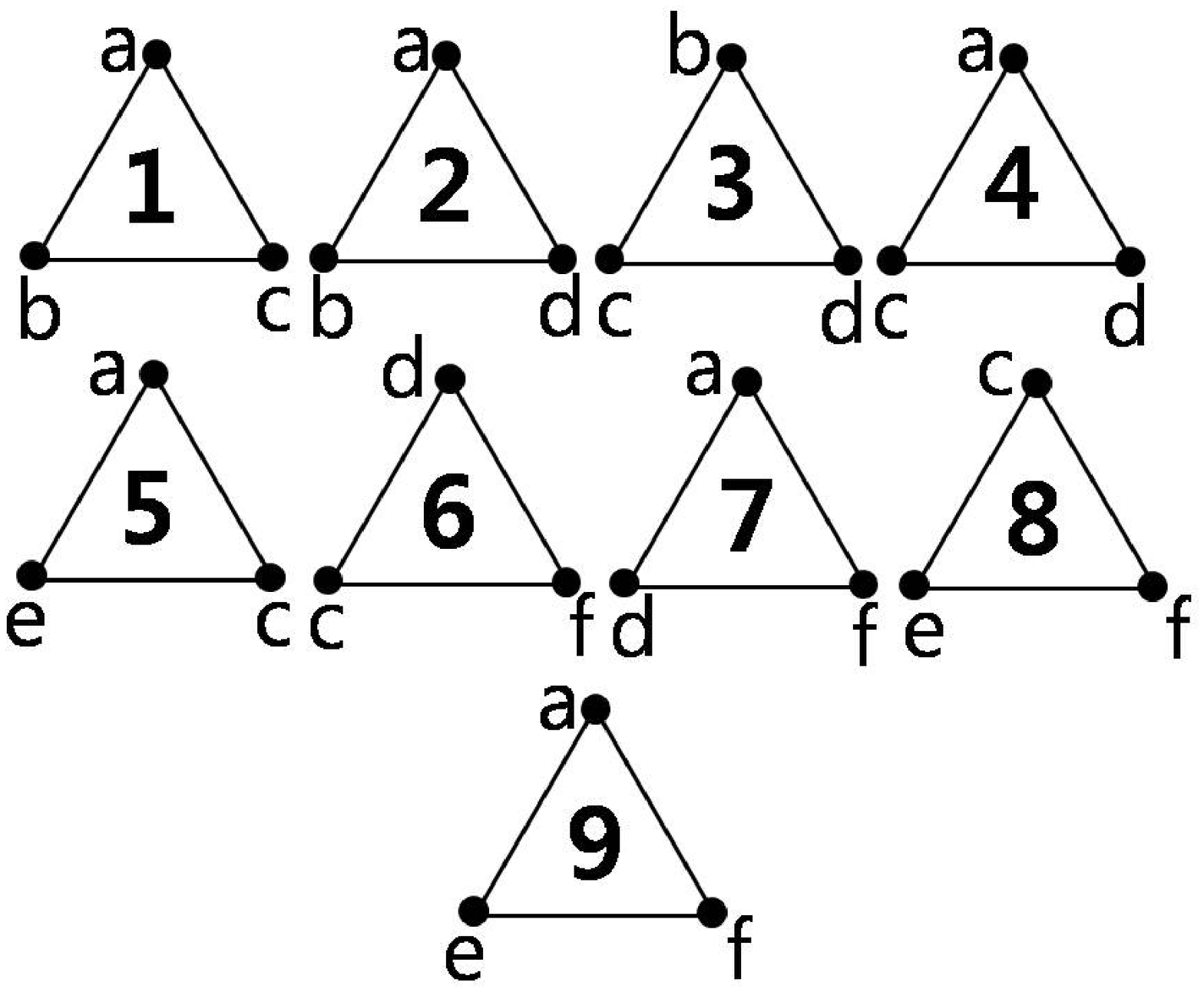,height=3cm,width=4cm}&\epsfig{file=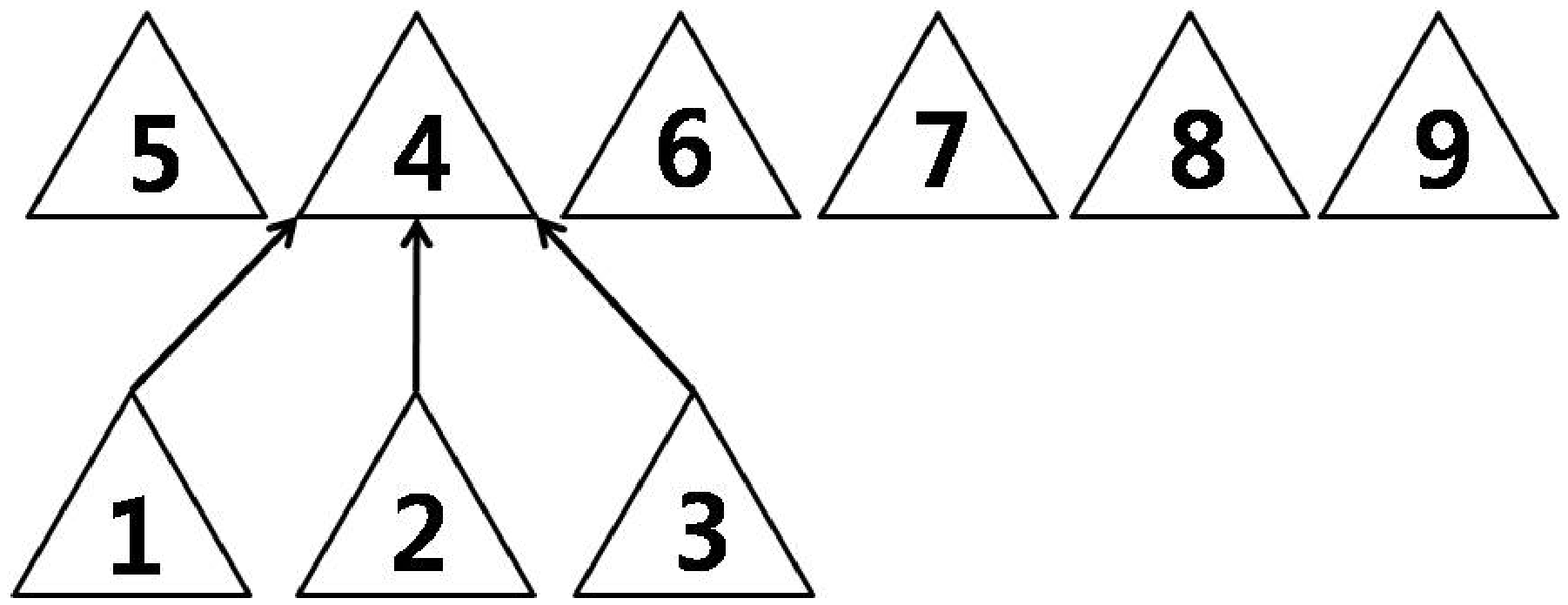,height=3cm,width=4cm}\\
\hbox{(a)}&\hbox{(b)}&\hbox{(c)}\\
\epsfig{file=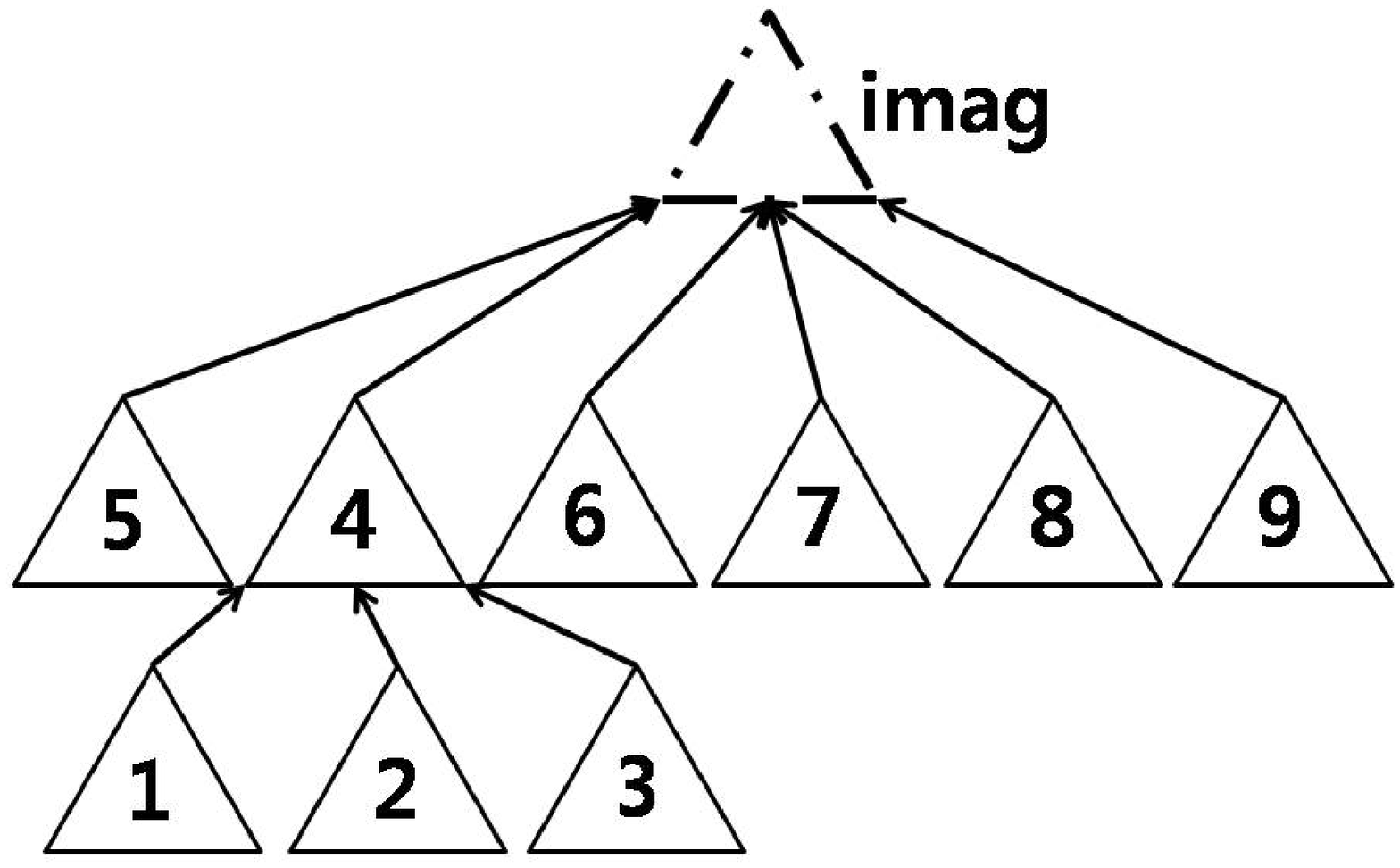,height=3cm,width=4cm}&\epsfig{file=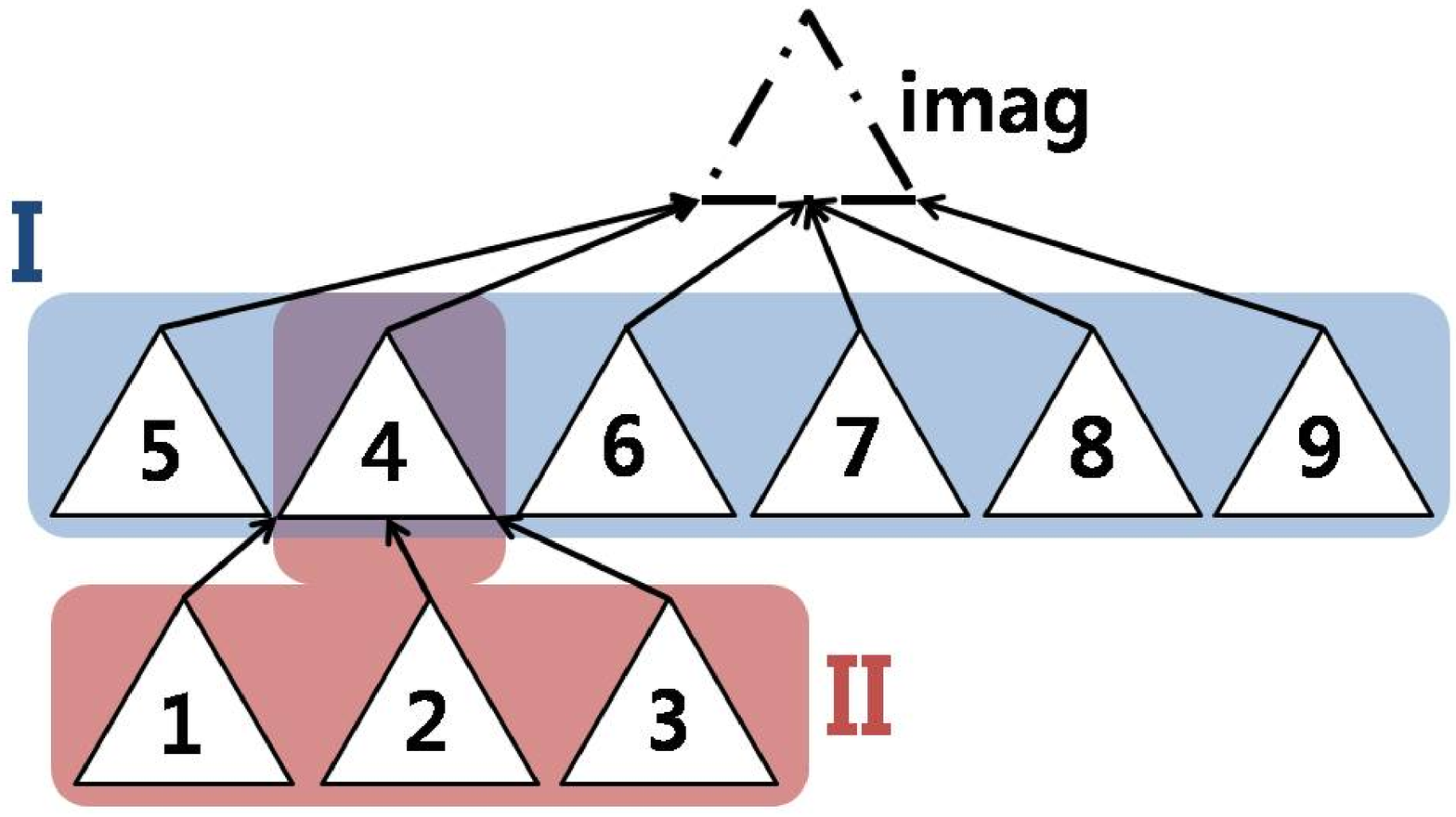,height=3cm,width=4cm}&\epsfig{file=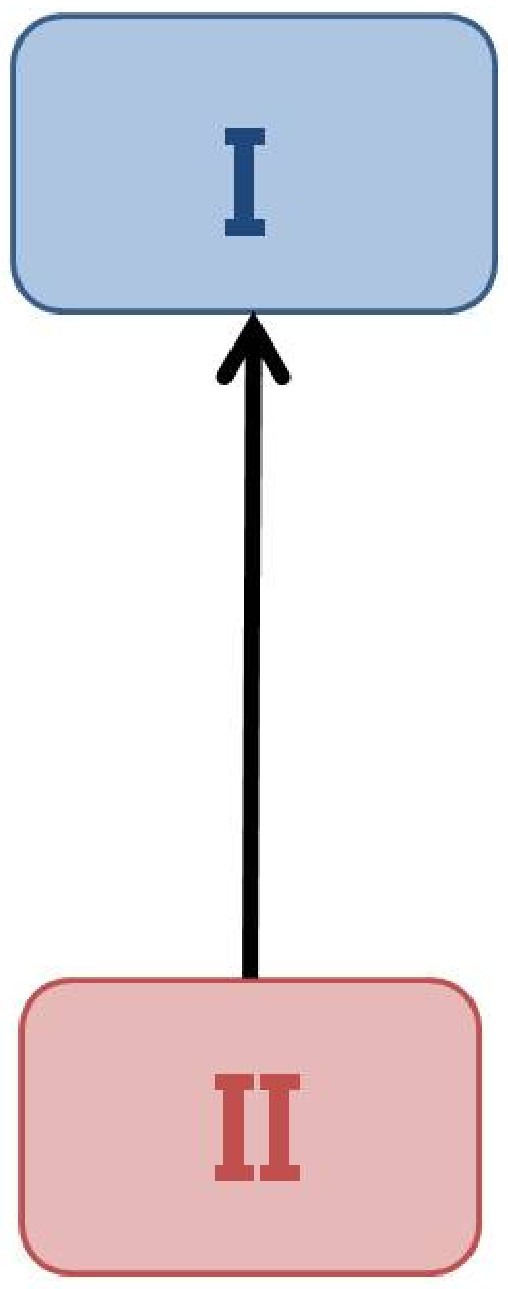,height=3cm,width=4cm}\\
\hbox{(d)}&\hbox{(e)}&\hbox{(f)}
\end{tabular}
\caption{An example of extraction of 3-clique and bubble hierarchies from a simple maximal planar graph.
(a) The maximal planar graph of interest.
(b) All 3-cliques.
(c) The 3-clique hierarchical tree.
(d) Insertion of $k_{imag}$ to merge the maximal 3-cliques producing the poset $(K',\preceq_{K'})$ and the corresponding hierarchical tree.
(e) 3-cliques which belong to different bubbles are highlighted.
(f) Bubble hierarchical tree. \label{SmallEg}}
\end{figure}

\section{Examples}\label{S.4}
In this section, we present two examples which will help to illustrate the relation between the graph structure and its hierarchical trees.
\subsection{Combination of Polyhedral Graphs}
In Fig.~\ref{SmallEg}, we have drawn a maximal planar graph which is made of six vertices (a,b,c,d,e and f). We first count and remunerate all 3-cliques in the graph as reported in Fig~\ref{SmallEg}(b). Then we can assign the relation between 3-cliques $\preceq_K$ by comparing the interiors as in Definition~\ref{Orient}. In this case, we can see that the 3-clique $(a,d,c)$ remunerated as `4' is the only 3-clique with a non-empty interior. Therefore, `4' has incoming neighbors in the hierarchical tree as depicted in Fig~\ref{SmallEg}(c).\\
We can now extend this to the bubble hierarchy. In Fig.~\ref{SmallEg}(e), the imaginary 3-clique (Definition~\ref{ImagCliq}) is included  so that the poset $(K',\preceq_{K'})$ has a single maximal element. By applying Theorems~\ref{BubMake} and~\ref{MaxBub}, we merge the 3-cliques with incoming neighbors in the hierarchical tree to obtain the list of bubbles. In  Fig.~\ref{SmallEg}(e), it is shown that we obtain 2 bubbles, namely  I and II. This results in the bubble hierarchical tree in Fig.~\ref{SmallEg}(f) according to the poset in Definition~\ref{BubHier}.\\
\begin{figure}[t]
\centering
\epsfig{file=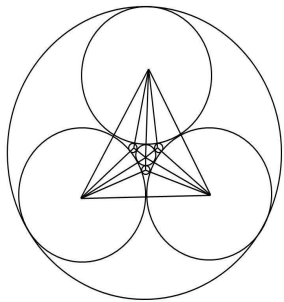,height=4cm,width=5cm}\caption{Apollonian packing at 3rd generation and the corresponding graph.\label{ApolloCircle}}
\end{figure}
\begin{figure}[]
\centering
\begin{tabular}{cc}
\epsfig{file=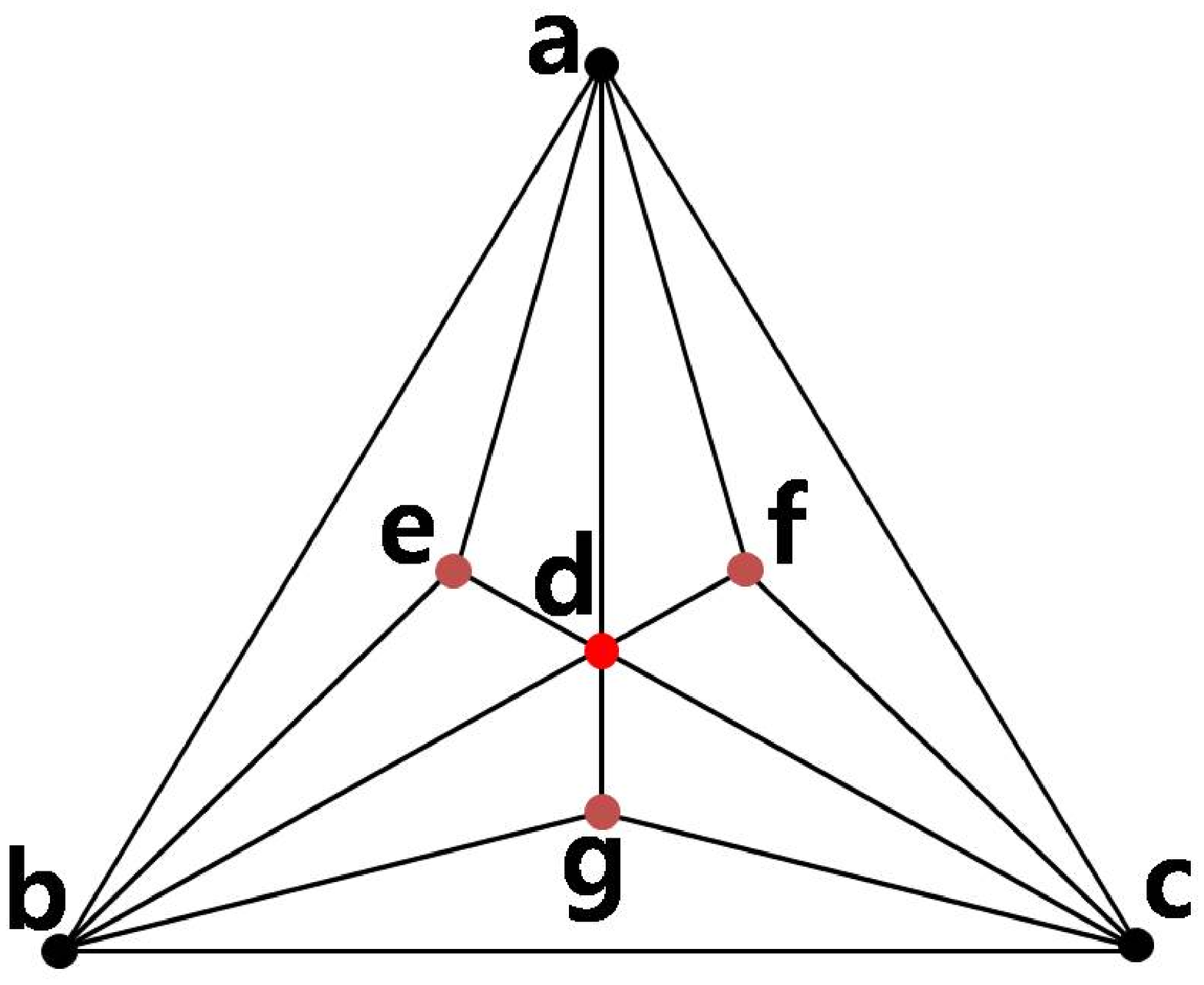,height=4cm,width=5cm}&\epsfig{file=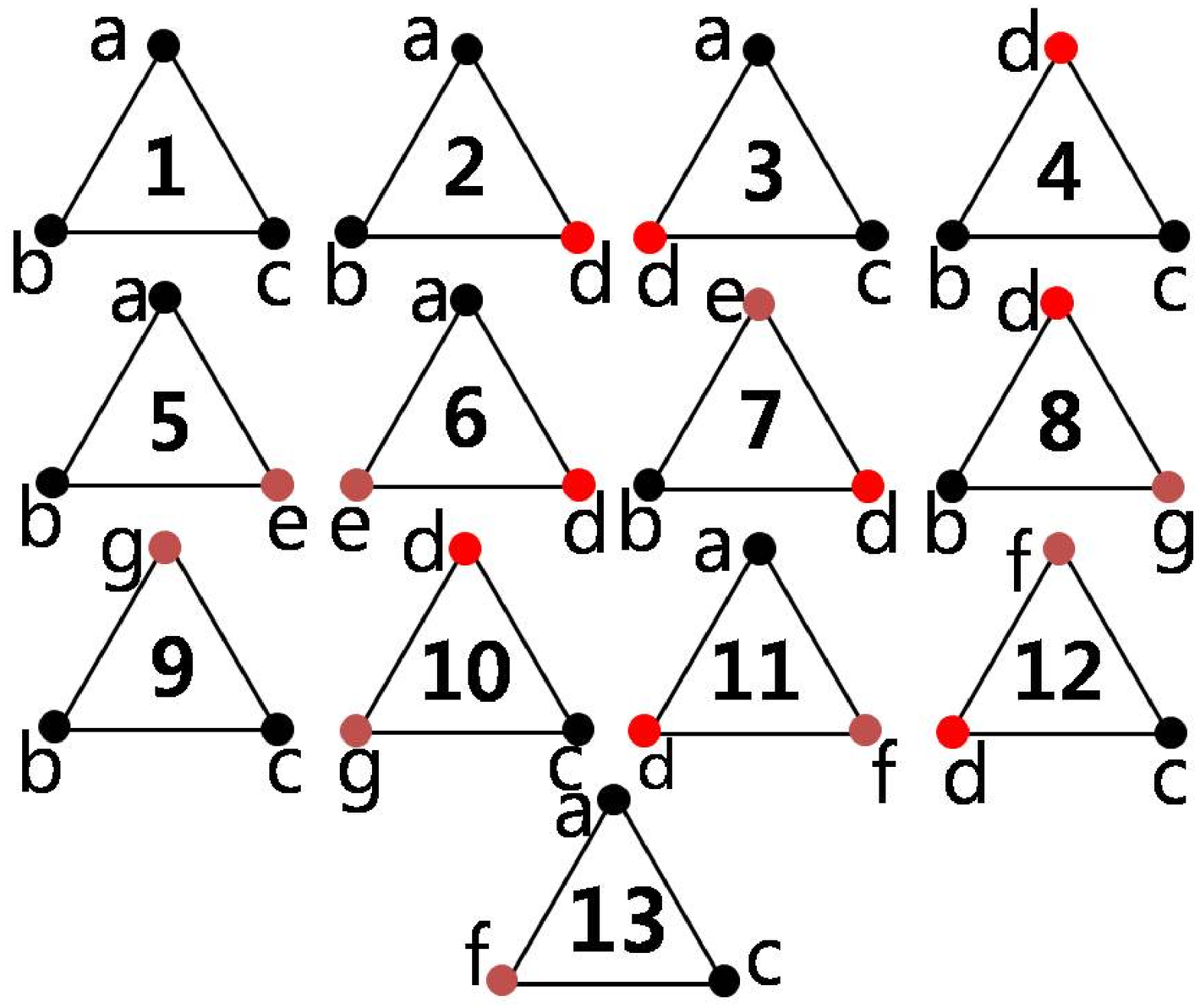,height=4cm,width=5cm}\\
\hbox{(a)}&\hbox{(b)}\\
\epsfig{file=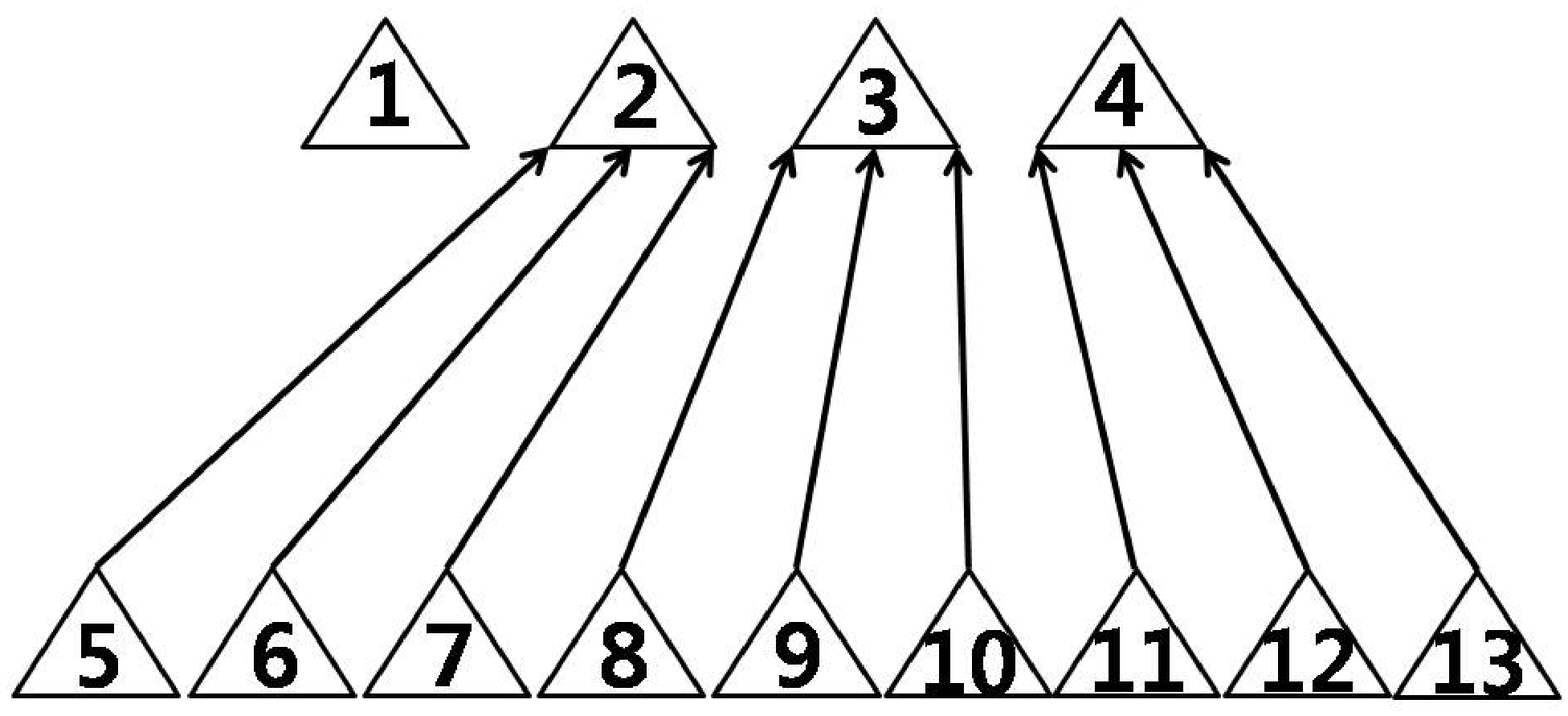,height=4cm,width=5cm}&\epsfig{file=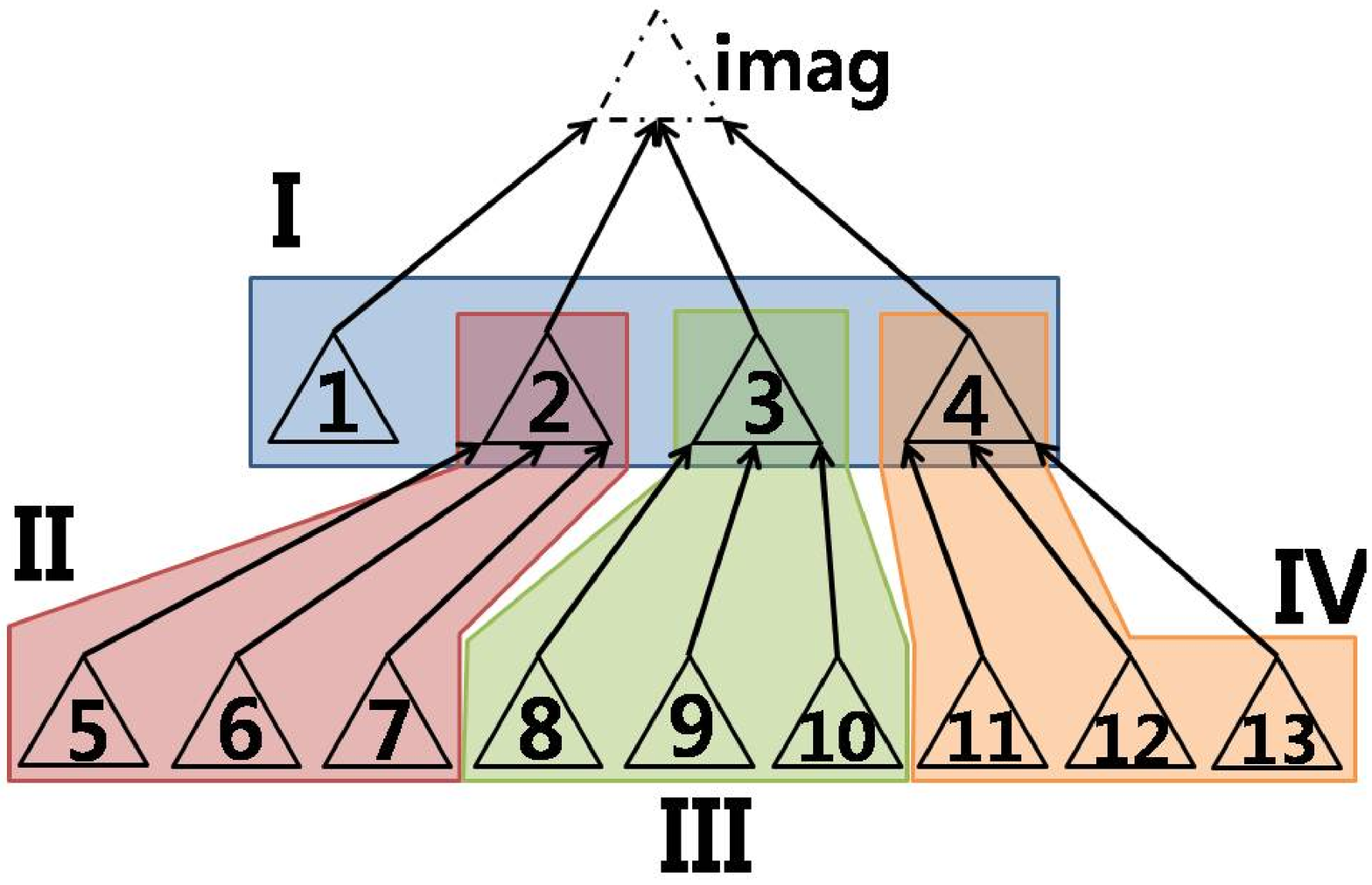,height=4cm,width=5cm}\\
\hbox{(c)}&\hbox{(d)}\\
\epsfig{file=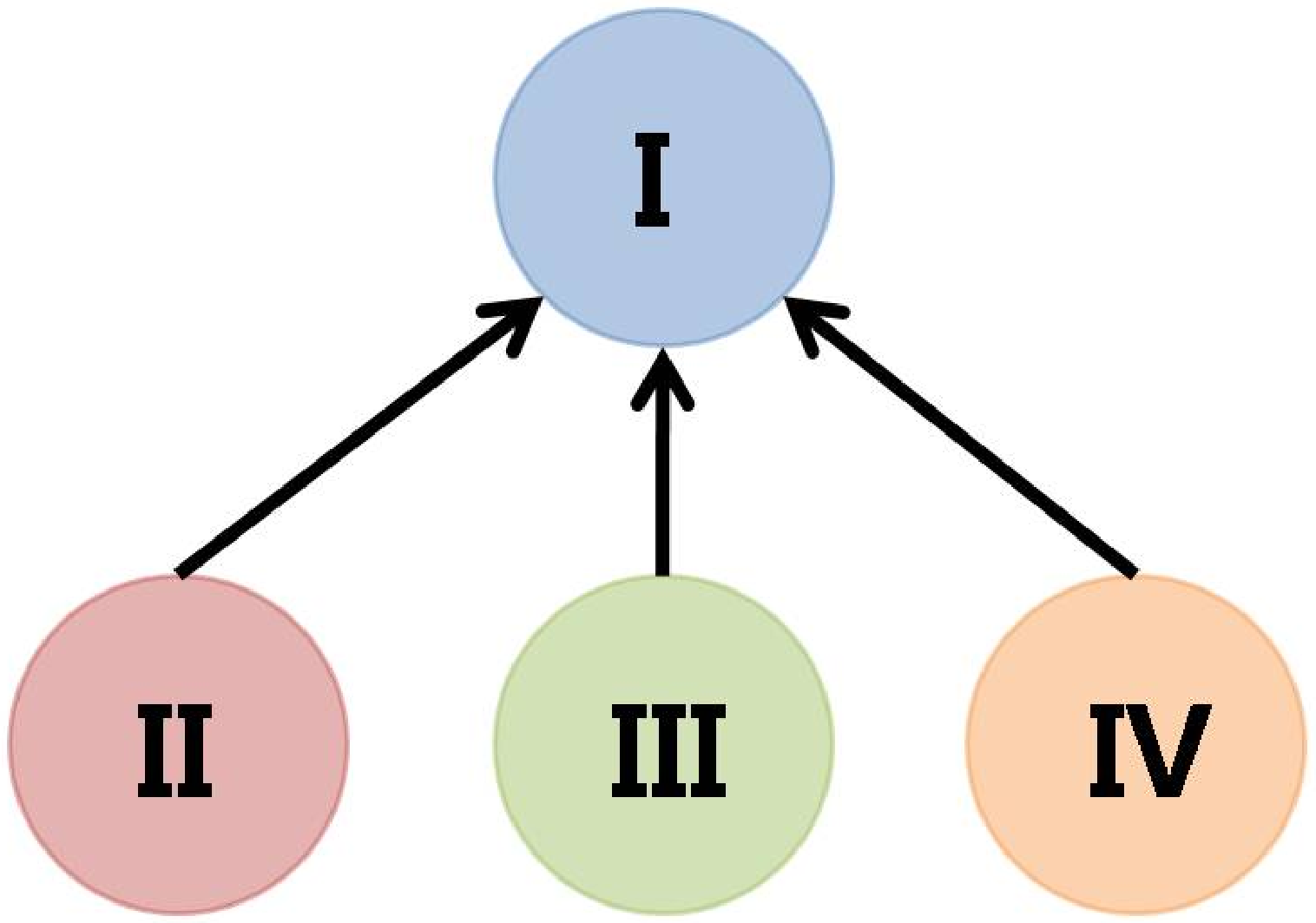,height=4cm,width=5cm}&\\
\hbox{(e)}&\\
\end{tabular}
\caption{(a) Apollonian Graph at $3$rd generation (b) all the 3-cliques. (c) 3-clique hierarchical tree.
(d) 3-clique hierarchical tree with the imaginary 3-clique.
(e) The bubble hierarchical tree.\label{Apollo}}
\end{figure}

\subsection{ Apollonian Graph }
In Fig.\ref{ApolloCircle} we report a type of maximal planar graph which has been inspired from the Apollonian packing of circles. An apollonian graph is constructed by connecting centers of tangent circles which are inserted at different stages to fill the voids \cite{AsteWeaire,Andrade2005}.
In particular, an Apollonian graph at $n$th generation is obtained from the Apollonian packing up to the insertion of the $n$th smallest circles. Beginning from a tetrahedral graph, this is equivalent to adding a vertex in each non-separating 3-clique, then join the extra vertex with the vertices of the 3-clique. In Fig.~\ref{ApolloCircle}, we have an example where the 3rd generation is reached and the corresponding graph is drawn.\\
On this graph (redrawn with labels in Fig.~\ref{Apollo}(a)) we have computed the corresponding hierarchical trees and bubbles as shown Fig.~\ref{Apollo}(b)-(e). From the figures, we see that all 3-cliques of nesting depth 0 have incoming neighbors except for `1' in (c).
Indeed, this Apollonian packing omits packing circles on the outside\footnote{With `outside' we refer to the exterior of the triangle `1' which has no vertices.} of the triangle `1' and does not allow it to have any incoming 3-cliques in $H$.
Similarly to Fig.~\ref{SmallEg}, we have classified the 3-cliques into the bubbles in (d) by making use of Theorems~\ref{BubMake} and \ref{MaxBub}. The bubble hierarchical tree is shown in (e).\\

\section{Conclusions}\label{S.5}
In this paper, we have shown that it is possible to define a unique hierarchical structure $(K,\preceq_K)$ of 3-cliques in maximal planar graphs $G$.
From the poset relation $\preceq_K$, we can build a hierarchical graph $H(K,E_k)$ which is a forest of rooted trees where the vertices are the 3-cliques $k_i\in K$ and there is a directed edge $\overrightarrow{k_jk_i}\in E_k$ from $k_j$ to $k_i$ if $k_i$ covers $k_j$ in $(K,\preceq_K)$.
This hierarchy depends on the definition of interior/exterior of the 3-clique.\\
We have shown that the extension of the 3-clique hierarchy to bubbles yields to a unique hierarchical structure for the bubbles. The set of bubbles $B$ from $G$ can be always uniquely identified and they form a tree $H_b(B',E_b)$ where the vertices are the bubbles $b_i \in B'$ and there is a directed edge $\overrightarrow{b_jb_i}\in E_b$ from $b_j$ to $b_i$ if $b_i$ covers $b_j$ in $(B',\preceq_{B'})$. We have shown that two neighboring bubbles are joined by one and only one 3-clique of $G$. The undirected topological structure of  $H_b(B,E_b)$ is independent from the definition of interior/exterior for the 3-cliques.\\
In the language of current network theory \cite{Newman2003}, it is natural to associate these bubbles with the idea of `communities'.
Indeed, they are connected portions of the graph which are loosely connected with the rest of the graph through 3-cliques.
The mathematical framework developed in this paper is therefore the base for a new way of identifying communities and detecting their relationships.
Applications to the analysis of weighted graphs and correlation based networks \cite{Tumminello2005,Aste2005,DiMatteo2005,Tumminello2007,DiMatteo2009} are under current investigation.

\section*{Acknowledgments}
This work was partially supported by the ARC Discovery Projects DP0344004 (2003), DP0558183 (2005) and COST MP0801 project.

\begin{appendix}
\section{Proof of  Theorem~\ref{BubMake}: Bubbles in $G$}\label{A.0}\ \\
\begin{proof}
\begin{enumerate}
\item[Connectedness]
   {Given a 3-clique $k_i$ with incoming neighbors $k_j,k_l,...$ in $H$, let us prove connectedness by proceeding in four steps:
    \begin{remunerate}
    \item define a set of disjoint subgraphs $\Omega=\{G^i_{out},G^j_{in},G^l_{in},...\}$;
    \item pick two subgraphs from $\Omega$, $G^j_{in}$ and $G^l_{in}$ and show that two vertices $v_j\in G^j_{in}$, $v_l\in G^l_{in}$ are connected in $G$ if the graph union $k_i\cup k_j\cup k_l\cup...$ is connected;
    \item prove that $(k_i\cup k_j\cup k_l\cup...)=G\setminus(G^i_{out}\cup G^j_{in}\cup G^l_{in}\cup...)$;
    \item prove connectedness of $k_i\cup k_j\cup k_l\cup...$.
    \end{remunerate}}
    {\bf First}, by Lemma~\ref{cliqueremoval}, $k_i$ distinguishes two disjoint subgraphs of $G$, $G^i_{in}$ and $G^i_{out}$ by removing $k_i$. This implies $G^i_{out}\cap(k_i\cup G^i_{in})=\emptyset$. Since $k_j,k_l,...$ are incoming neighbors of $k_i$, $(k_j\cup G^j_{in})\subseteq(k_i\cup G^i_{in})$. Therefore, $G^j_{in}\cap G^i_{out}=\emptyset$ for any incoming neighbor $k_j$ at $k_i$. Moreover, by Theorem~\ref{SeparateCliq}, $G^j_{in}\cap G^l_{in}=\emptyset$ for any two incoming neighbors $k_j$ and $k_l$. Therefore, the set of subgraphs $\Omega$ consists of disjoint graphs.\\
    {\bf Second}, since $k_l\not\in G^j_{in}$ and $k_j\not\in G^l_{in}$, $k_j$ and $k_l$ must be connected in $G\setminus(G^j_{in}\cup G^l_{in})$ in order to maintain connectedness of $G$. Repeating this argument for all the pairs of subgraphs in $\Omega$, this yields that $k_i,k_j,k_l,...$ are connected in $G\setminus(G^i_{out}\cup G^j_{in}\cup G^l_{in}\cup...)$.\\
{    {\bf Third}, we now must prove that $(k_i\cup k_j\cup k_l\cup...)=G\setminus(G^i_{out}\cup G^j_{in}\cup G^l_{in}\cup...)$.
Let us suppose that there is a vertex $v$ such that $v\not\in(k_i\cup k_j\cup k_l\cup...)$ and $v\in G\setminus(G^i_{out}\cap G^j_{in}\cap G^l_{in}\cap...)$ so that we prove the claim if existence of such $v$ yields contradiction.
Since $G$ is maximally planar, $v$ belongs to at least one 3-clique. Let us call this 3-clique $k_m$. We know that $k_m\neq k_j$ for any incoming neighbor $k_j$ by $v\not\in(k_i\cup k_j\cup k_l\cup...)$. By Corollary~\ref{InCliq}, $k_m\subseteq(k_i\cup G^i_{in})$, hence $k_m\preceq_K k_i$. On the other hand, $k_m\not\preceq_K k_j$ for any incoming neighbor $k_j$ because $v\not\in G^j_{in}\Rightarrow k_m\not\subseteq(k_j\cup G^j_{in})\Rightarrow (G^j_{in}\cap G^m_{in})=\emptyset$.
Therefore, $k_m\preceq_K k_i$ and $k_m\not\preceq_K k_j$.
Consequently, $k_i$ is the only 3-clique that satisfies $k_m\preceq_K k_i$, hence $k_i$ covers $k_m$. Then $k_m$ is another incoming neighbor at $k_i$ in $H$ by the Definition~\ref{H}. However, $k_j,k_l,...$ are all of the incoming neighbors as we have assumed from the beginning. This yields a contradiction, therefore the claim is true.}\\
{\bf Fourth} having proven that $(k_i\cup k_j\cup k_l\cup...)=G\setminus(G^i_{out}\cup G^j_{in}\cup G^l_{in}\cup...)$, we are now in a position to say that the graph union $b_i=(k_i\cup k_j\cup k_l\cup...)$ is connected since we have already proved $G\setminus(G^i_{out}\cup G^j_{in}\cup G^l_{in}\cup...)$ is connected.

\item[Maximally Planar]
Firstly, let us stress that we have already proved that $G$ is a disjoint union $G=b_i\cup(G^i_{out}\cup G^j_{in}\cup G^l_{in}\cup...)$. This implies that if one could add an extra edge in $b_i$ without violating the planarity, one can also do it in $b_i\cup(G^i_{out}\cup G^j_{in}\cup G^l_{in}\cup...)$ since this is a disjoint union. This is not possible because $G$ is maximally planar.
Therefore $b_i$ is maximally planar as well.
\item[Non-Separating] Since $b_i$ involves only the 3-cliques $k_i$, $k_j$, $k_l,...$, and does not include $G^i_{out},G^j_{in},G^l_{in},...,$ this implies that each 3-clique in $b_i$ has either its interior of exterior removed from $G$. Therefore each 3-clique does not separate $b_i$ into two subgraphs since either interior or exterior is always empty. Therefore they are non-separating 3-cliques.
\end{enumerate}
\end{proof}

\section{Proof of Maximal Bubble in Theorem~\ref{MaxBub}}\label{A.1}
Here we must prove that the graph union of maximal 3-cliques of $(K,\preceq_K)$ is a bubble and its root 3-clique is $k_{imag}$ (Theorem ~\ref{MaxBub}).
We do this by following the same reasoning used in the proof of Theorem~\ref{BubMake}.
\begin{proof}
\begin{enumerate}
\item[Connectedness]$b_{\rho}$ must be connected. Suppose it is not so that there are disconnected components of $b_{\rho}$. That is, there exists at least one pair of vertices in $v_p,v_q\in V_{\rho}$ which are not connected. Let us pick an arbitrary vertex $v_p$ from one of the disconnected components and another $v_q$ from a different component. This implies $v_p$ and $v_q$ belong to different maximal 3-cliques $k_{\rho_i}$ and $k_{\rho_j}$. Because $k_{\rho_i}$ and $k_{\rho_j}$ are separating 3-cliques, any paths between $v_p$ and $v_q$ must contain at least one vertex from each of $k_{\rho_i}$ and $k_{\rho_j}$. In order to maintain the connectedness of $G$, $k_{\rho_i}$ and $k_{\rho_j}$ must be connected.
\item[Maximally Planar] We claim that $b_{\rho}=(\bigcup_i k_{\rho_i})=(G\setminus(\bigcup_i G^{\rho_i}_{in}))$. The argument is very similar to that of $b_i=(k_i\cup k_j\cup k_l\cup...)=G\setminus(G^i_{out}\cup G^j_{in}\cup G^l_{in}\cup...)$ in the proof of Theorem~\ref{BubMake} in Appendix~\ref{A.0}. Suppose there exists $v\not\in(\bigcup_i k_{\rho_i})$ but $v\in(G\setminus(\bigcup_i G^{\rho_i}_{in}))$. Since $G$ is maximally planar, there exists $k_i$ such that $v\in k_i$ and $k_i\neq k_{\rho_i}$ for all $\rho_i$.\\
Since $k_{\rho_i}$ are all of the maximal 3-cliques in $(K,\preceq_K)$,
$k_{\rho_i}\not\preceq_K k_i$ for all $\rho_i$, but $k_i\preceq_K k_{\rho_i}$ for some $\rho_i$. Then $k_i\subseteq(k_{\rho_i}\cup G^{\rho_i}_{in})\Rightarrow v\in G^{\rho_i}_{in}$. This is against the initial assumption that $v\in(G\setminus(\bigcup_i G^{\rho_i}_{in}))$. Therefore, $(\bigcup_i k_{\rho_i})=(G\setminus(\bigcup_i G^{\rho_i}_{in}))$. Then, one can say $G$ is a disjoint union of $b_{\rho}\cup(\bigcup_i G^{\rho_i}_{in})$.\\
Let us call $\overline{b_{\rho}}=(\bigcup_i G^{\rho_i}_{in})$ for simplicity. Then, any additional edge in $b_{\rho}$ without violating the planarity can be added in $G$ as well, since $G$ is the disjoint union $b_{\rho}\cup\overline{b_{\rho}}$. Since $G$ is maximally planar, this is not true. Therefore, no additional edge can be added in $b_{\rho}$ without violating the planarity.
\item[Non-Separating] Following the same argument in the proof of Theorem~\ref{BubMake}, $b_{\rho}=G\setminus\overline{b_{\rho}}$, hence all 3-cliques in $b_{\rho}$ do not have any interior. Therefore they are non-separating in $b_{\rho}$.
\end{enumerate}
\end{proof}
\section{Proof on the set of bubbles in Corollary \ref{AllBub}}\label{A.2}
Let $K_b=\{k_i,k_j,k_l,...\}$ be the set of 3-cliques of a given bubble $b$.
Before proceeding to prove Corollary~\ref{AllBub}, let us state some useful theorems.
\begin{theorem}[Removal of 3-clique from a bubble]\label{RemovalonBub}
Given a bubble $b$ and a 3-clique $k_i\subset b$, $(b\setminus k_i)\subseteq G^i_{in}$ or $(b\setminus k_i)\subseteq G^i_{out}$.
\end{theorem}
\begin{proof}
By Lemma~\ref{cliqueremoval}, there are two subgraphs $G^i_{in}$ and $G^i_{out}$ such that $(G\setminus k_i)=(G^i_{in}\cup G^i_{out})$ is a disjoint union. Then any connected subgraph of $(G\setminus k_i)$ must be a subgraph of either $G^i_{in}$ or $G^i_{out}$. Therefore, $(b\setminus k_i)\subseteq G^i_{in}$ or $(b\setminus k_i)\subseteq G^i_{out}$ since $(b\setminus k_i)$ is connected as $k_i$ is non-separating in $b$.
\end{proof}
The following corollary is immediate from Theorem~\ref{RemovalonBub}.
\begin{corollary}[Maximal 3-clique of $K_b$]\label{NumMax}
Given $k_i,k_j\in K_b$, there exist at most one 3-clique $k_i$ in $K_b$ such that, for any $k_j\in K_b$, $k_j\preceq_K k_i$. Moreover, $k_i$ covers $k_j$.
\end{corollary}
\begin{proof}
It is immediate from Theorem~\ref{RemovalonBub} that, if there exists $k_i$ such that $(b\setminus k_i)\subseteq G^i_{in}$, then all $k_j\preceq_K k_i$ for all $k_j\in K_b$.\\
Now, suppose that $k_i$ is the 3-clique with property $k_j\preceq_K k_i$ for all $k_j\in K_b$. Then, suppose there exists $k_j\in K_b$ that $k_i$ does not cover $k_j$. Then, there is another 3-clique $k_l$ such that $k_l\preceq_K k_i$ and $k_j\preceq_K k_l$, so $k_i\subseteq(k_l\cup G^l_{out})$ and $k_j\subseteq(k_l\cup G^l_{in})$. If $k_l\in K_b$, then this implies $k_l$ is a separating 3-clique in $b$, therefore $b$ is not a bubble. On the other hand, if $k_l\not\in K_b$, then $b$ is a disconnected subgraph, hence $b$ is not a bubble either. By contradiction, $k_i$ must cover $k_j$.
\end{proof}
Now, we extend the latter statement in Corollary~\ref{NumMax} that no single maximal 3-clique exists in $K_b$.
\begin{corollary}[Maximal bubble]\label{BubMax}
If there is no $k_i\in K_b$ such that $k_j\preceq_K k_i$ for all other $k_j\in K_b$, then $b$ is the bubble made of maximal 3-cliques in the poset $(K,\preceq_K)$.
\end{corollary}
\begin{proof}
Suppose that there is no $k_i\in K_b$ such that $k_j\preceq_K k_i$ for all other $k_j\in K_b$. And, suppose there exists $k_j\in K_b$ that $k_j$ is not maximal in $(K,\preceq_K)$. Then there is $k_l\in K$ that $k_j\preceq_K k_l$. If $k_l\in K_b$, then $b$ is not a bubble since $k_l$ is a separating 3-clique in $b$. On the other hand, if $k_l\not\in K_b$, then $b$ is a disconnected subgraph of $G$, therefore $b$ is not a bubble either. Therefore, by contradiction, there is no such $k_j\in K_b$ that is not maximal in $(K,\preceq_K)$.
\end{proof}
Now, let us state the final result of Theorem~\ref{RemovalonBub} and Corollaries~\ref{NumMax} and \ref{BubMax}.
\begin{corollary}[Bubbles in $H$]\label{BubAll}
If $b$ is a bubble in $G$, then the set of 3-cliques in $K_b$ are either:
\begin{remunerate}
\item[(i)]Union of a 3-clique $k_i$ and all of its incoming neighbors $k_j,k_l,...$ in $H$, or
\item[(ii)]Union of all maximal 3-cliques of $(K,\preceq_K)$.
\end{remunerate}
\end{corollary}
\begin{proof}
By Corollaries~\ref{NumMax} and \ref{BubMax}, a bubble $b$ consists of a set $K_b$ of 3-cliques that consists of (i) one maximal 3-clique $k_i$ and covered 3-cliques $k_j$ (i.e. incoming neighbors in $H$ at $k_i$), or (ii) maximal 3-cliques. What the corollaries do not prove is whether $K_b$ is the set of (i) $k_i$ and \textbf{all} of covered 3-cliques $k_j$, or (ii)\textbf{all} maximal 3-cliques.\\
Let us suppose that $K_b$ does not include all incoming neighbors at $k_i$. And let us call the set of 3-cliques $k_i$ and all of its incoming neighbors as $K_i$. We have shown that union of 3-cliques in $K_i$ yields a bubble in the proof of Theorem~\ref{BubMake}. Since $K_b\subset K_i$, this implies union of 3-cliques in $K_b$ is not maximally planar. Therefore, the assumption that $K_b\neq K_i$ is wrong. Therefore $K_b=K_i$.\\
Similarly, the union of all maximal 3-cliques yields a maximally planar bubble as stated in Theorem~\ref{MaxBub}. Therefore, Corollary~\ref{BubAll} is true.
\end{proof}
It is immediate that  Corollary~\ref{BubAll} is the equivalent statement to Corollary~\ref{AllBub}. Therefore we have proven Corollary~\ref{AllBub}.
\end{appendix}

\bibliographystyle{plain}

\begin{thebibliography}{}

\end{thebibliography}


\begin{thebibliography}{10}

\bibitem{Caldarelli2007}
G. Caldarelli,
``Scale-Free Networks. Complex Webs in Nature and Technology'',
Oxford University Press, 2007.

\bibitem{Boccaletti2006}
S.~Boccaletti, V.~Latora, Y.~Moreno, M.~Chavez, and D.-U. Hwang,
\newblock Complex networks: Structure and dynamics,
\newblock {\em Phys. Rep.}, 424:175--308, 2006.

\bibitem{Girvan2002}
M.~Girvan and M.~E.~J. Newman,
\newblock Community structure in social and biological networks,
\newblock {\em Proc. Natl. Acad. Sci. USA}, 99(12):7821--7826, 2002.

\bibitem{Aaron2008}
Aaron Clauset, Cristopher Moore, and M.E.J. Newman,
\newblock Hierarchical structure and the prediction of missing links in
  networks,
\newblock {\em Nature}, 453:98--101, 2008.

\bibitem{Newman2003}
M.~Girvan M.E.J.~Newman,
\newblock Finding and evaluating community structure in networks,
\newblock {\em Phys. Rev. E}, 69(026113), 2004.

\bibitem{Tumminello2005}
Michel Tumminello, Tomaso Aste, Tiziana~Di Matteo, and Rosario.~N. Mantegna,
\newblock A tool for filtering information in complex systems,
\newblock {\em Proc. Natl. Acad. Sci. USA}, 102(30), 2005.

\bibitem{Aste2005}
T~Aste, T.~Di Matteo, and S~Hyde,
\newblock Complex networks on hyperbolic surfaces,
\newblock {\em Phys. A}, 346:20--26, 2005.

\bibitem{DiMatteo2005}
T.~Di Matteo, T.~Aste, S.~T. Hyde, and S.~Ramsden,
\newblock Interest rates hierarchical structure,
\newblock {\em Phys. A}, 21--33, 2005.

\bibitem{Tumminello2007}
M.~Tumminello, T.~Di Matteo, T.~Aste, and R.N. Mantegna,
\newblock Correlation based networks of equity returns sampled at different
  time horizons.
\newblock {\em Eur. Phys. J. B}, 55(2):209--217, 2007.

\bibitem{DiMatteo2009}
T.~Di~Matteo, F.~Pozzi, T.~Aste,
\newblock The use of dynamical networks to detect the hierarchical organization of financial market sectors,
\newblock {\em Eur. Phys. J. B} 2009.

\bibitem{Diestel2005}
Reinhard Diestel,
\newblock {\em Graph Theory ed.3},
\newblock Springer-Verlag, 2005.

\bibitem{Andrade2005}
Jos\'{e} S.~Andrade Jr., Hans~J. Herrmann, Roberto~F.S. Andrade, and Luciano~R.
  da~Silva,
\newblock Apollonian networks: Simultaneously scale-free, small world,
  euclidan, space filling and matching graphs,
\newblock {\em Phys. Rev. Lett.}, 94(018702), 2005.

\bibitem{Jech2003}
Jech. T.,
\newblock {\em Set Theory},
\newblock Springer-Verlag, 2003.

\bibitem{AsteWeaire}
T.~Aste and D.~Weaire,
\newblock {\em The Pursuit of Perfect Packing},
\newblock Taylor and Francis, 2nd Ed., New York, 2008.

\end{thebibliography}
\end{document}